\documentclass[11pt,a4paper]{article}

\usepackage[utf8]{inputenc}
\usepackage[english]{babel}
\usepackage{amsmath}
\usepackage{amsthm}
\usepackage{amsfonts}
\usepackage{amssymb}
\usepackage{amscd}
\usepackage{tikz-cd}
\usepackage{bezier}
\usepackage{enumerate}
\usepackage{fancyhdr}
\usepackage{graphicx}
\usepackage{epstopdf}
\usepackage{booktabs}
\usepackage{rotating}
\usepackage{listings}
\usepackage[square, numbers, comma, sort&compress]{natbib} 
\usepackage{fancyhdr}
\usepackage{amsmath,amsfonts,amssymb,amscd,amsthm,xspace}
\usepackage[centerlast,small,sc]{caption}
\usepackage[pdfpagemode={UseOutlines},bookmarks=true,bookmarksopen=true,
   bookmarksopenlevel=0,bookmarksnumbered=true,hypertexnames=false,
   colorlinks,linkcolor={blue},citecolor={blue},urlcolor={red},
   pdfstartview={FitV},unicode,breaklinks=true]{hyperref}

\newcommand{\bigslant}[2]{{\raisebox{.2em}{$#1$}\left/\raisebox{-.2em}{$#2$}\right.}}
\newtheorem{lemma}{Lemma}
\newtheorem{proposition}{Proposition}
\newtheorem*{proposition*}{Proposition}
\newtheorem*{theorem*}{Theorem}
\newtheorem{theorem}{Theorem}

\newtheorem{definition}{Definition}

\title{Cylindrical contact homology and topological entropy }
\author{Marcelo Alves \\ \textsl{Universit\'e Libre de Bruxelles \& Universit\'{e} Paris-Sud}}

\begin{document}
\maketitle
\begin{center}
\abstract{\textbf{Abstract.} We establish a relation between the growth of the cylindrical contact homology of a contact manifold and  the topological entropy of Reeb flows on this manifold. We show that if a contact manifold $(M,\xi)$ admits a hypertight contact form $\lambda_0$  for which the cylindrical contact homology has exponential homotopical growth rate, then the Reeb flow of every contact form on $(M,\xi)$ has positive topological entropy. Using this result, we provide numerous new examples of contact 3-manifolds on which every Reeb flow has positive topological entropy.}
\end{center}

\section{Introduction} \label{introduction}
The aim of this paper is to establish a relation between the behaviour of cylindrical contact homology and the topological entropy of Reeb flows. The topological entropy is a non-negative number associated to a dynamical system which measures the complexity of the orbit structure of the system. Positivity of the topological entropy means that the system possesses some type of exponential instability. We show that if the cylindrical contact homology of a contact 3-manifold is ``complicated enough'' from a homotopical viewpoint, then every Reeb flow on this contact manifold has positive topological entropy.

\subsection{Basic definitions and history of the problem}

We first recall some basic definitions from contact geometry. A 1-form $\lambda$  on a $(2n+1)$-dimensional manifold $Y$ called a \textit{contact form} if $ \lambda \wedge (d\lambda)^n $ is a volume form on $Y$. The hyperplane $\xi= \ker \lambda$ is called the \textit{contact structure}. For us a \textit{contact manifold} will be a pair $(Y,\xi)$ such that $\xi$ is the kernel of some contact form $\lambda$ on $Y$ (these are usually called co-oriented contact manifolds in the literature). When $\lambda$ satisfies $\xi = \ker\lambda$, we will say that $\lambda$ is a contact form on $(Y,\xi)$. On any contact manifold there always exist infinitely many different contact forms. Given a contact form $\lambda$, its
\textit{Reeb vector field} is the unique vector field $X_\lambda$ satisfying $\lambda(X_\lambda)=1$ and $i_{X_\lambda}d\lambda=0$. The Reeb flow of $\lambda$ is the flow of the vector field $X_\lambda$. We will refer to the periodic orbits of the Reeb flow as \textit{Reeb orbits}.

We study the topological entropy of Reeb flows from the point of view of contact topology. More precisely, we search for conditions on the topology of a contact manifold $(M,\xi)$ that force \textbf{all} Reeb flows on $(M,\xi)$ to have positive topological entropy. The condition we impose is on the behaviour of a contact topological invariant called cylindrical contact homology. We show that if a contact manifold $(M,\xi)$ admits a contact form $\lambda_0$ for which the cylindrical contact homology has \textit{exponential homotopical growth}, then all Reeb flows on  $(M,\xi)$ have positive topological entropy.

The notion of exponential homotopical growth of cylindrical contact homology, which is introduced in this paper, differs from the notion of growth of contact homology studied in \cite{CH,Vaugon}. For reasons explained in Section~\ref{section2}, the growth of contact homology is not well adapted to study the topological entropy of Reeb flows, while the notion of homotopical growth rate is (as we show) well suited for this purpose.
We begin by explaining the results which were previously known relating the behaviour of contact topological invariants to the topological entropy of Reeb flow.

The study of contact manifolds all of whose Reeb flows have positive topological entropy was initiated by Macarini and Schlenk \cite{MS}. They showed that if $Q$ is an energy hyperbolic manifold and $\xi_{geo}$ is the contact structure on the unit tangent bundle $T_1 Q$ associated to the geodesic flows, then every Reeb flow on $(T_1 Q,\xi_{geo})$ has positive topological entropy. Their work was based on previous ideas of Frauenfelder and Schlenk \cite{FS1,FS2} which related the growth rate of Lagrangian Floer homology to entropy invariants of symplectomorphisms. The strategy to estimate the topological used in \cite{MS} can be briefly sketched as follows:
\begin{center}
Exponential growth of Lagrangian Floer homology of the tangent fiber $(TQ)_{|_p}$ \\
$\Rightarrow$ \\
Exponential volume growth of the unit tangent fiber $(T_1 Q)_{|_p}$ for all Reeb flows in $(T_1 Q,\xi_{geo})$ \\
$\Rightarrow$ \\
Positivity of the topological entropy for all Reeb flows in $(T_1 Q,\xi_{geo})$.
\end{center}
To obtain the first implication Macarini and Schlenk use the fact that $(T_1 Q,\xi_{geo})$ has the structure of a Legendrian fibration, and apply the geometric idea of \cite{FS1,FS2} to show that the number of trajectories connecting a  Legendrian fiber to  another Legendrian fiber can be used to obtain a volume growth estimate.
The second implication in this scheme follows from Yomdin's theorem, that states that exponential volume growth of a submanifold implies positivity of topological entropy. \footnote{The same scheme was used in \cite{FS3,FLS} to obtain positive lower bounds for the intermediate and slow entropies of Reeb flows on unit tangent bundles; we discuss these results in more detail in Section~\ref{section7}.}

In the author's Ph.D. thesis \cite{A3,A1} this approach was extended to deal with $3$-dimensional contact manifolds which are not unit tangent bundles. This was done by designing a localised version of the geometric idea of \cite{FS1,FS2}. Globally most contact $3$-manifolds are not Legendrian fibrations, but a small neighbourhood of a given Legendrian knot in any contact $3$-manifold can be given the structure of a Legendrian fibration. It turns out that this is enough to conclude that if the linearised Legendrian contact homology of a pair of Legendrian knots in a contact 3-manifold $(M^3,\xi)$  grows exponentially, then the length of these Legendrian knots grows exponentially for any Reeb flow on $(M^3,\xi)$. We then apply Yomdin's theorem to obtain that all Reeb flows on $(M^3,\xi)$ have positive topological entropy.

One drawback of these approaches is that they only give lower entropy bounds for $C^{\infty}$-smooth Reeb flows. The reason  is that Yomdin's theorem holds only for $C^{\infty}$-smooth flows. The approach presented in the present paper \textbf{does not} use Yomdin's theorem and gives lower bounds for the topological entropy of $C^1$-smooth Reeb flows.

Another advantage is that the cylindrical contact homology is usually easier to compute than the linearised Legendrian contact homology. In fact, to apply the strategy of \cite{A3,A1} to a contact 3-manifold $(M^3,\xi)$ one must first find a pair of Legendrian curves which, one believes, ``should'' have exponential growth of linearised Legendrian contact homology. This is highly non-trivial since on any contact 3-manifolds there exist many Legendrian links for which the linearised Legendrian contact homology does not even exist. On the other hand the definition of cylindrical contact homology only involves the contact manifold $(M^3,\xi)$, and no Legendrian submanifolds.

\subsection{Main results}

Our results are inspired by the philosophy that a ``complicated'' topological structure can force chaotic behavior for dynamical systems associated to this structure. Two important examples of this phenomena are: the fact that on manifolds with complicated loop space the geodesic flow always has positive topological entropy  (see \cite{Pat}), and the fact that every diffeomorphism of a surface which is isotopic to a pseudo-Anosov diffeomorphism has positive topological entropy \cite{Fel}.

To state our results we introduce some notation. Let $M$ be a manifold and $X$ be a $C^k$ ($k \geq 1$ ) vector field. Our first result relates the topological entropy of $\phi_X$ to the growth (relative to $T$) of the number of distinct homotopy classes which contain periodic orbits of $\phi_X$ with period $\leq T$. More precisely let  $\Lambda_X^T$ be the set of free homotopy classes of $M$ which contain a periodic orbit of $\phi_X$ with period $\leq T$. We denote by $N_X(T)$ the cardinality of  $\Lambda_X^T$.

\begin{theorem} \label{theorem1}
If for real numbers $a > 0$ and $b$ we have $N_X(T) \geq  e^{aT+b}$, then $h_{top}( X)\geq a$.
\end{theorem}

Theorem~\ref{theorem1} might be a folklore result in the theory of dynamical systems. However as we have not found it in the literature, we provide a complete proof in Section~\ref{section2}. It contains as a special case Ivanov's inequality for surface diffeomorphisms (see \cite{J}).
Our motivation for proving this result is to apply it to Reeb flows. Contact homology allows one to carry over information about the dynamical behaviour of one special Reeb flow on a contact manifold to all other Reeb flows  on the same contact manifold. In Section~\ref{section4} we introduce the notion of exponential homotopical growth of cylindrical contact homology. As we already mentioned, this growth rate differs from the ones  previously considered in the literature and is specially designed to allow one to use Theorem~\ref{theorem1} to obtain results about the topological entropy of Reeb flows. This is made via the following:

\begin{theorem} \label{theorem2}
Let $\lambda_0$ be a hypertight contact form on a contact manifold $(M,\xi)$ and assume that the cylindrical contact homology $C\mathbb{H}_{cyl}^{J_0}(\lambda_0)$ has exponential homotopical growth with exponential weight $a>0$. Then for every $C^k$ ($k\geq2$) contact form $\lambda$ on $(M,\xi)$ the Reeb flow of $X_\lambda$ has positive topological entropy. More precisely, if $f_\lambda$ is the unique function such that $\lambda = f_\lambda \lambda_0$, then
\begin{equation}
h_{top}(X_{\lambda})\geq \frac{a}{\max f_\lambda}.
 \end{equation}
\end{theorem}
Notice that Theorem~\ref{theorem2} allows us to conclude the positivity of the topological entropy for \textbf{all} Reeb flows on a given contact manifold $(M,\xi)$, once we show that  $(M,\xi)$ admits one special hypertight contact form for which the cylindrical contact homology has exponential homotopical growth. It is worth remarking that our proof of Theorem~\ref{theorem2} is carried out in full rigor, and does \textbf{not} make use of the Polyfold technology which is being  developed by Hofer, Wysocki and Zehnder. The reason is that we do not use the linearised contact homology  considered in \cite{BEE,Vaugon}, but resort to a topological idea used in \cite{HMS} to prove existence of Reeb orbits in prescribed homotopy classes.

Theorem~\ref{theorem2} above allows one to obtain estimates for the topological entropy for $C^1$-smooth Reeb flows. As previously observed, the strategy used in \cite{MS,A3,A1} produces estimates for the topological entropy only for $C^{\infty}$-smooth contact forms as they depend on Yomdim's theorem, which fails for finite regularity.

Our other results are concerned with the existence of examples of contact manifolds which have a contact form with exponential homotopical growth rate of cylindrical contact homology. We show that in dimension 3 they exist in abundance, and it follows from Theorem \ref{theorem2} that every Reeb flow on these contact manifolds has positive topological entropy. In Section~\ref{section5} we construct such examples for manifolds with a non-trivial JSJ decomposition and with a hyperbolic component that fibers over the circle.

\begin{theorem} \label{theorem3}
Let $M$ be a closed connected oriented  3-manifold which can be cut along a nonempty family of incompressible tori into a family $\{M_i, 0 \leq i \leq k\}$ of irreducible manifolds with boundary, such that the component $M_0$ satisfies:
\begin{itemize}
\item{$M_0$ is the mapping torus of a diffeomorphism $h: S \to S$ with pseudo-Anosov monodromy on a surface $S$ with non-empty boundary.}
\end{itemize}
Then $M$ can be given infinitely many non-diffeomorphic contact structures $\xi_k$, such that for each $\xi_k$ there exists a hypertight contact form $\lambda_k$ on $(M,\xi_k)$ which has exponential homotopical growth of cylindrical contact homology.
\end{theorem}

In Section~\ref{section6} we study the cylindrical contact homology of contact 3-manifolds $(M,\xi_{(q,\mathfrak{r})})$ obtained via a special integral Dehn surgery on the unit tangent bundle $(T_1 S,\xi_{geo})$ of a hyperbolic surface $(S,g)$. This Dehn surgery is performed on a neighbourhood of a Legendrian curve $L_{\mathfrak{r}}$ which is the Legendrian lift of a separating geodesic. The surgery we consider is the contact version of Handel-Thurston surgery  which was introduced by Foulon and Hasselblatt in \cite{FH} to produce non-algebraic Anosov Reeb flows in 3-manifolds. We call this contact surgery the Foulon-Hasselblatt surgery. This surgery produces not only a contact 3-manifold $(M,\xi_{(q,\mathfrak{r})})$, but also a special contact form which we denote by $\lambda_{FH}$ on $(M,\xi_{(q,\mathfrak{r})})$.
 In \cite{FH} the authors restrict their attention to integer surgeries with positive surgery coefficient $q$ and prove that, in this case, the Reeb flow of $\lambda_{FH}$ is Anosov. Our methods work also for negative coefficients as the Anosov condition on $\lambda_{FH}$ does not play a role in our results. We obtain:

\begin{theorem} \label{theorem4}
Let $(M,\xi_{(q,\mathfrak{r})})$ be the contact manifold endowed obtained by the Foulon-Hasselblatt surgery, and $\lambda_{FH}$ be the contact form obtained via the Foulon-Hasselblat surgery on the Legendrian lift $L_{\mathfrak{r}} \subset T_1 S$. Then $\lambda_{FH}$ is hypertight and its cylindrical contact homology has exponential homotopical growth.
\end{theorem}

\textbf{Organization of the paper.} In Section~\ref{section2} we recall one of the definitions of the topological entropy and present the proof of Theorem \ref{theorem1}. In Section~\ref{section3} we recall the definition of cylindrical contact homology and its basic properties. In Section~\ref{section4} we introduce the notion of exponential homotopical growth of cylindrical contact homology and prove Theorem \ref{theorem2}. Section~\ref{section5} is devoted to the proof of Theorem \ref{theorem3}. In Section~\ref{section6} we present the definition of the integral Foulon-Hasselblatt surgery and prove Theorem \ref{theorem4}. In Section~\ref{section7} we discuss the results obtained in this paper and propose some questions for future research.

\textbf{Remark:} We again would like to point out that all the results above \textbf{do not} depend on the Polyfolds technology which is being developed Hofer, Wysocki and Zehnder. This is the case because the versions of contact homology used  for proving the results above involve only somewhere injective pseudoholomorphic curves. In this situation transversality can be achieved by ``classical'' perturbation methods as in \cite{Dr}.

\

\textbf{Acknowledgements:} I specially thank my professors Fr\'ed\'eric Bourgeois and Chris Wendl for their guidance, support and for our many discussions which were crucial for the development of this paper, which is a part of my PhD thesis being developed under their supervision. I would like to thank professor Pedro Salom\~ao for many helpful discussions and for explaining to me the techniques used in \cite{HMS} which made it possible to avoid dealing with transversality problems arising from multiply covered pseudoholomorphic curves. My thanks to professor Felix Schlenk for his interest in this work, his suggestions for improving the exposition and for suggesting many directions for future work. My personal thanks to Ana Nechita, Andr\'e Alves, Hilda Ribeiro and Lucio Alves for their unconditional personal support. Lastly, I thank FNRS-Belgium for the financial support.

\setcounter{theorem}{0}

\section{Homotopic growth of periodic orbits and topological entropy} \label{section2}

Throughout this section $M$ will denote a compact manifold. We endow $M$ with an auxiliary Riemannian metric $g$, which induces a distance function $d_g$ on $M$, whose injective radius we denote by $\epsilon_g$. Let $\widetilde{M}$ be the universal cover of $M$, $\widetilde{g}$ be the Riemannian metric that makes the covering map $\pi: \widetilde{M} \to M$ an isometry, and $d_{\widetilde{g}}$ be the distance induced by the metric $\widetilde{g}$.

Let $X$ be a vector field on $M$ with no singularities and $\phi_{X}^t$ the flow generated by $X$. We call $P^{X}(T)$ the number of periodic orbits of $\phi^t$ with period in $[0,T]$. For us a periodic orbit of $X$ is a pair $([\gamma]_c,T)$ where $[\gamma]_c$ is the set of parametrizations of a given \emph{immersed} curve $c:S^1 \to M$, and $T$ is a positive real number (called the period of the orbit) such that:

\begin{itemize}
\item{$\gamma \in [\gamma]_c \ \iff $ $\gamma: \mathbb{R} \to M$ parametrizes $c$ and $\dot{\gamma}(t)= X(\gamma(t))$ }
\item{for all $\gamma \in [\gamma]_c$ we have $\gamma(T+t)= \gamma(t)$ and $\gamma([0,T])= c$ }
\end{itemize}
We say that a periodic orbit $([\gamma]_c,T)$ is in a free homotopy class $l$ of $M$ if $c \in l$.

By a parametrized periodic orbit $(\gamma,T)$ we mean a periodic orbit $([\gamma]_c,T)$ with a fixed choice of parametrization $\gamma \in [\gamma]_c$. A parametrized periodic orbit $(\gamma,T)$ is said to be in a free homotopy class $l$ when the underlying periodic orbit $([\gamma]_c,T)$ is in $l$.

We now recall a definition of topological entropy due to Bowen \cite{Bo} which will be very useful for us. Let $T$ and $\delta$ be positive real numbers. A set $S$ is said to be $T,\delta$-separated if for all $q_1 \neq q_2 \in S$ we have:
\begin{equation}
\max_{t\in [0,T]}d_{g}(\phi_X^t(q_1),\phi_X^t(q_2))> \delta.
\end{equation}
We denote by $n^{T,\delta}$ the maximal cardinality of a $T,\delta$-separated set for the flow $\phi_X$. Then we define the $\delta$-entropy $h_\delta(\phi_X)$ as:
\begin{equation}
h_{\delta}(\phi_X)=\limsup_{T\to +\infty } \frac{\log(n^{T,\delta})}{T}
\end{equation}
The topological entropy $h_{top}$ is then defined by
\begin{equation}
h_{top}(\phi_X)= \lim_{\delta \to 0} h_{\delta}(\phi_X).
\end{equation}
One can prove that the topological entropy does not depend on the metric $d_g$ but only on the topology determined by the metric. For these and other structural results about topological entropy we refer the reader to any standard textbook in dynamics such as \cite{HK} and \cite{Rob}.

From the work of Kaloshin and others it is well known that the exponential growth rate of periodic orbits $\limsup_{T \to +\infty}\frac{\log(P^{X}(T))}{T}$ can be much bigger than the topological entropy. This implies that the growth rate $\limsup_{T \to +\infty}\frac{\log(P^{X}(T))}{T}$ does not give a lower bound for the topological entropy of an arbitrary flow. There is however a different growth rate, which measures how quickly periodic orbits appear in different free homotopy classes, and which can be used to give such a lower bound of the topological entropy of a flow.

Let $\Lambda$ denote the set of free homotopy classes of loops in $M$, and $\Lambda_0 \subset \Lambda$ the subset of primitive free homotopy classes. Denote by $\Lambda^T_X \subset \Lambda$ the set of free homotopy classes $\varrho$ such that there exists a periodic orbit of $\phi^t_X$ with period smaller or equal to $T$ which is homotopic to $\varrho$. We denote by $N_{X}(T)$ the cardinality of $\Lambda^T_X$.

Let $\{(\gamma_i, T_i) ; 1 \leq i \leq n\}$ be a finite set of parametrized periodic orbits of $X$. For a number $T$ satisfying $T\geq T_i$ for all $i \in \{1,...,n\}$ and a constant $\delta >0$, we denote by $\Lambda^{T,\delta}_X((\gamma_1,T_1),..., (\gamma_n,T_n)) $ the subset of $\Lambda$ such that:

\begin{itemize}
\item{$l \in \Lambda^{T,\delta}_X((\gamma_1,T_1),..., (\gamma_n,T_n))$ if, and only if, there exist a parametrized periodic orbit $(\widehat{\gamma},\widehat{T})$ with period $\widehat{T} \leq T$ in the free homotopy class $l$ and a number $i_l \in \{1,...,n\}$ for which $\max_{t\in [0,T]} (d_g(\gamma_{i_l}(t),\widehat{\gamma}(t) )) \leq \delta$ . }
\end{itemize}

Notice that
\begin{equation}
\Lambda^{T,\delta}_X((\gamma_1,T_1),..., (\gamma_n,T_n))= \bigcup_{i \in \{1,...,n\}} \Lambda^{T,\delta}_X((\gamma_i,T_i)).
\end{equation}

We are ready to prove the main result in this section. Theorem \ref{theorem1'} below is well known to be true in the particular cases where $\phi_X$ is a geodesic flow, where it follows from Manning's inequality (see \cite{K2} and \cite{Pat}); and where $\phi_X$ is the suspension of surface diffeomorphism with pseudo-Anosov monodromy, where it follows from Ivanov's theorem (see \cite{J}). It can be seen as a generalization of these results in the sense that it includes them as particular cases and that it applies to many other situations. Our argument is inspired by the remarkable proof of Ivanov's inequality given by Jiang in (\cite{J}).

\begin{theorem} \label{theorem1'}
If for real numbers $a>0$ and $b $ we have $N_{X}(T)\geq e^{aT+b}$, then $h_{top}(\phi_X)\geq a$.
\end{theorem}

\textit{Proof:}
The theorem will follow if we prove that for all $\delta < \frac{\epsilon_g}{10^6}$ we have $h_{\delta}(\phi_X)\geq a$. From now on fix $0<\delta< \frac{\epsilon_g}{10^6}$.

\textbf{Step 1:}  For any point $p \in M$ let $V_{4\delta}(p)$ be the $4\delta$-neighbourhood of $\pi^{-1}(p)$. Because $\delta < \frac{\epsilon_g}{10^6}$, it is clear that $V_{4\delta}(p)$ is the disjoint union

\begin{equation}
V_{4\delta}(p)=\bigcup_{\widetilde{p}\in \pi^{-1}(p) } B_{4\delta}(\widetilde{p} )
\end{equation}
where the ball $B_{4\delta}(\widetilde{p})$ is taken with respect to the metric $\widetilde{g}$.

Because of our choice of $\delta< \frac{\epsilon_g}{10^6}$ it is clear that there exists a constant $0<k_1$ which does not depend on $p$, such that if $B$ and $B'$ are two distinct connected components of $V_{4\delta}(p)$ we have $d_{\widetilde{g}}(B,B')> k_1$.

Because of compactness of $M$, we know that the vector field $\widetilde{X}:= \pi^* X$ is bounded in the norm given by the metric $\widetilde{g}$. Combining this with the inequality in the last paragraph, one obtains the existence of a constant $0<k_2$, which again doesn't depend on $p$ such that, if $\widetilde{\upsilon}:[0,R] \to \widetilde{M} $ is a parametrized trajectory of $\phi_{\widetilde{X}}$ such that $\widetilde{\upsilon}(0) \in B$ and $\widetilde{\upsilon}(R) \in B'$ for $B\neq B'$ are connected components of $V_{4\delta}(p)$ then $R> k_2$.

From the last assertion we deduce the existence of a constant $\widetilde{K}$, depending only $g$ and $X$, such that for every $p\in M$ and every parametrized trajectory $\widetilde{\upsilon}:[0,T] \to \widetilde{M} $ of $\phi_{\widetilde{X}}$, the number $L^T(p,\widetilde{\upsilon})$ of distinct connected components of $V_{4\delta}(p)$, intersected by the curve $\widetilde{\upsilon}([0,T]) $  satisfies

\begin{equation}
L^T(p,\widetilde{\upsilon}) < \widetilde{K} T + 1.
\end{equation}

\textbf{Step 2:} We claim that for every parametrized periodic orbit $(\gamma',T')$ of $X$ we have
\begin{equation}
 \sharp(\Lambda^{T,\delta}_X((\gamma',T')) < \widetilde{K}T + 1
\end{equation}
for all $T>T'$.

To see this take $\widetilde{\gamma}'$ be a lift of $\gamma'$ and let $p'=\gamma'(0)$ and $\widetilde{p}' = \widetilde{\gamma'}(0)$.
We consider the set $\{B_j; 1 \leq j \leq \mathfrak{m}^T(\gamma',T')\}$ of connected components of $V_{4\delta}(p')$ satisfying:

\begin{itemize}
\item{$B_j \neq B_k$ if $j \neq k$,}
\item{if $B$ is a connected component of $V_{4\delta}(p')$ which intersects $\widetilde{\gamma'}([0,T])$ then $B=B_j$ for some $j \in \{1,..., \mathfrak{m}^T(\gamma',T')\} $,}
\item{if $j < i$ then $B_j$ is visited by the trajectory $\widetilde{\gamma'}: [0,T] \to  \widetilde{M}$ before $B_i$.}
\end{itemize}
From step 1, we know that $\mathfrak{m}^T(\gamma',T')< \widetilde{K} T +1$.

\begin{center}
\begin{minipage}[b]{0.9\textwidth}
    \includegraphics[width=\textwidth]{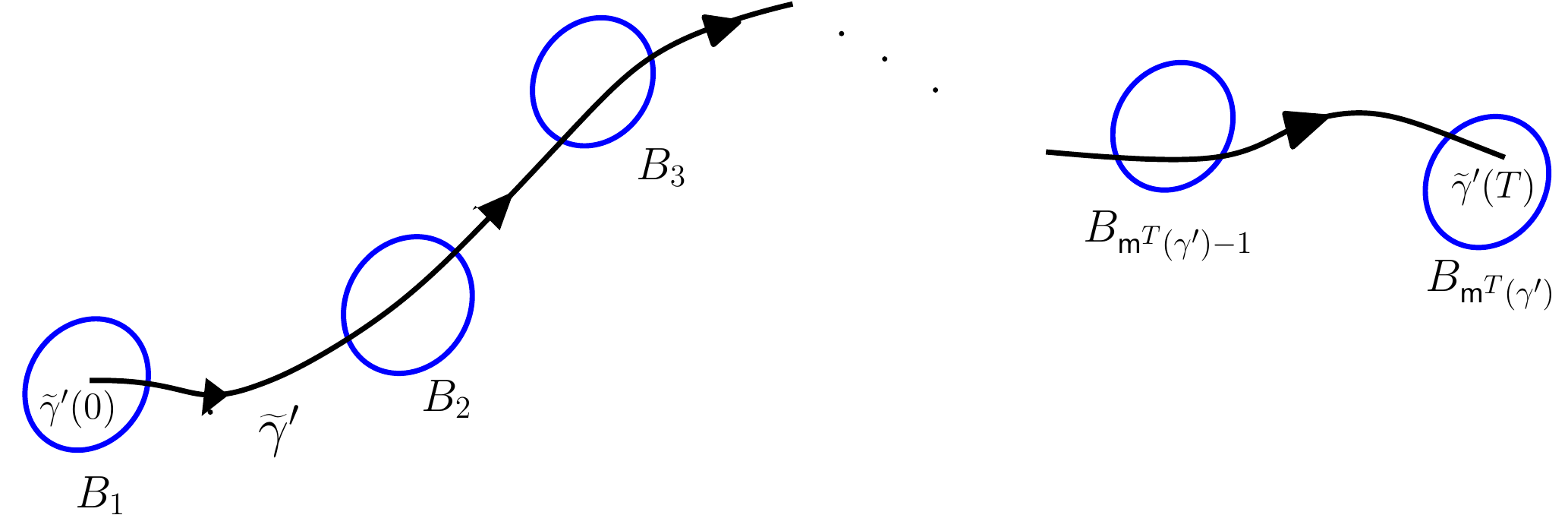}
    \label{fig:2}
    \begin{center}
      Figure 1: The set $\{B_j; 1 \leq j \leq \mathfrak{m}^T(\gamma',T')\}$.
    \end{center}
\end{minipage}
\end{center}

For each $l \in \Lambda^{T,\delta}_X((\gamma',T'))$ pick $(\chi_l,T_l)$ in $l$ to be a parametrized periodic orbit which satisfies $d_{g}(\chi_l(t),\gamma'(t))<\delta$ for all $t\in [0,T]$. There exists a lift $\widetilde{\chi_l}$ of $\chi_l$ satisfying $d_{\widetilde{g}}(\widetilde{\chi_l}(t),\widetilde{\gamma'}(t))<\delta$ for all $t\in [0,T]$.

From the triangle inequality it is clear that the point $q_l = \widetilde{\chi_l}(0)$ is in the connected component $B_1$ which contains $\widetilde{p}'$.
We will show that $\widetilde{\chi_l}(T_l)$ is contained in $B_j$ for some $j \in \{1,..., \mathfrak{m}^T(\gamma')\} $. Because $\pi(\widetilde{\chi_l}(0)) = \pi(\widetilde{\chi_l}(T_l))$, we have:

\begin{equation}
d_{\widetilde{g}}(\widetilde{\chi_l}(T_l),\pi^{-1}(p'))= d_{\widetilde{g}}(\widetilde{\chi_l}(0),\pi^{-1}(p')) < \delta
\end{equation}
which already implies that $\widetilde{\chi_l}(T_l) \in V_{4\delta}(p')$. We denote by $\widetilde{p'_l}$ the unique element $\pi^{-1}(p')$ for which we have $d_{\widetilde{g}}(\widetilde{\chi_l}(T_l),\widetilde{p_l}') < \delta$. Using the triangle inequality we now obtain:
\begin{equation}
d_{\widetilde{g}}(\widetilde{\gamma'}(T_l),\widetilde{p'_l}) \leq d_{\widetilde{g}}(\widetilde{\gamma'}(T_l),\widetilde{\chi_l}(T_l)) + d_{\widetilde{g}}(\widetilde{\chi_l}(T_l),\widetilde{p_l}') < \delta + \delta.
\end{equation}
From the inequalities above we conclude that $\widetilde{\gamma'}(T_l)$ and $\widetilde{\chi_l}(T_l)$ are in the connected component of $V_{4\delta}(p')$ that contains $\widetilde{p_l}'$. Because this connected component contains $\widetilde{\gamma'}(T_l)$, it is therefore one of the $B_j$ for $j \in \{1,..., \mathfrak{m}^T(\gamma',T')\}$ as we wanted to show. We can thus define a map $\Upsilon^{T,\delta}_{(\gamma',T')}: \Lambda^{T,\delta}_X((\gamma',T')) \to \{1,..., \mathfrak{m}^T((\gamma',T'))\}$ which associates to each $l \in \Lambda^{T,\delta}_X(\gamma')$ the unique $j \in \{1,..., \mathfrak{m}^T(\gamma',T')\}$ for which $\widetilde{\chi_l}(T_l) \in B_j$.

We now claim that if $l \neq l'$ then $\widetilde{\chi_l}(T_l)$ and $\widetilde{\chi_{l'}}(T_{l'})$ are in different connected components of $V_{4\delta}(p')$. To see this notice that both $\widetilde{\chi_l}(0)$ and $\widetilde{\chi_{l'}}(0)$ are in the component $B_1$. Therefore it is clear, because $\delta < \frac{\epsilon_g}{10^6}$, that if $\widetilde{\chi_l}(T_l)$ and $\widetilde{\chi_{l'}}(T_{l'})$ are in the same component of $V_{4\delta}(p')$, then the closed curves $\chi_l ([0,T_l])$ and $\chi_{l'}([0,T_{l'}])$ are freely homotopic. This contradicts our choice of $(\chi_l,T_l)$ and $(\chi_{l'},T_{l'})$ and the fact that $l\neq l'$.

We thus conclude that the map $\Upsilon^{T,\delta}_{(\gamma',T')}: \Lambda^{T,\delta}_X((\gamma',T')) \to \{1,..., \mathfrak{m}^T(\gamma',T')\}$ is injective, which implies that  $\sharp(\Lambda^{T,\delta}_X((\gamma',T'))) \leq \mathfrak{m}^T(\gamma',T') < \widetilde{K}T + 1 $.

\

\textbf{Step 3:} Inductive step.

As an immediate consequence of step 2 we have that if $\{(\gamma_i,T_i); 1 \leq i \leq m\} $ is a set of parametrized periodic orbits of $X$ we have $\sharp(\Lambda^{T,\delta}_X((\gamma_1,T_1),..., (\gamma_m,T_m))) \leq m(\widetilde{K}T+1)$.

\textbf{Inductive claim:} Fix $T>0$ and suppose that $S^{T}_{m}=\{(\gamma_i,T_i) ; 1 \leq i \leq m  \}$ is a set of parametrized periodic orbits such that $T\geq T_i$ for every $i\in \{1,...,m\}$, and that satisfies:

\begin{itemize}
\item{(a) The free homotopy classes $l_i$ of $(\gamma_i,T_i)$ and $l_j$ of $(\gamma_j,T_j)$ are distinct if $i\neq j$,}
\item{(b) For every $i\neq j$ we have $\max_{t\in [0,T]}d_{g}(\gamma_i(t),\gamma_j(t))> \delta$.}
\end{itemize}
Then, if $m < \frac{N_X(T)}{\widetilde{K}T + 1}$, there exists a parametrized periodic orbit $(\gamma_{m+1},T_{m+1}\leq T)$ such that its homotopy class $l_{m+1}$ does not belong to the set $\{l_i ; 1\leq i \leq m\}$ and such that
\begin{equation}
\max_{t\in [0,T]}d_{g}(\gamma_{m+1}(t),\gamma_i(t))> \delta
\end{equation}
for all $i\in {1,...,m}$.

\textit{Proof of the claim:} Recall that $\sharp(\Lambda^{T,\delta}_X((\gamma_1,T_1),..., (\gamma_m,T_m))) \leq m(\widetilde{K}T +1)$. Therefore, because $m < \frac{N_X(T)}{\widetilde{K}T +1 }$, there exists a free homotopy $l_{m+1} \in \Lambda^T_X \setminus \Lambda^{T,\delta}_X((\gamma_1,T_1),..., (\gamma_m,T_m))$. Choose a parametrized periodic orbit $(\gamma_{m+1},T_{m+1})$ with $T_{m+1}\leq T$ in the homotopy class $l_{m+1}$.

As $l_{m+1} \notin \Lambda^{T,\delta}_X((\gamma_1,T_1),..., (\gamma_m,T_m))$, we must have $\max_{t\in [0,T]}d_{g}(\gamma_{m+1}(t),\gamma_i(t))> \delta$ for all $i\in {1,...,m}$; thus completing the proof of the claim.

\

\textbf{Step 4:} Obtaining a $T,\delta$ separated set.

As usual, we denote by $ \lfloor \frac{N_X(T)}{\widetilde{K}T + 1} \rfloor$ the largest integer which is  $\leq  \frac{N_X(T)}{\widetilde{K}T + 1} $.
The strategy is now to use the inductive step to obtain a set $S^T_X =\{(\gamma_i,T_i); 1 \leq i \leq \lfloor \frac{N_X(T)}{\widetilde{K}T + 1} \rfloor \} $ satisfying conditions (a) and (b) above with the maximum possible cardinality. We start with a set $S^T_1 = \{(\gamma_1,T_1)\}$, which clearly satisfies  conditions (a) and (b), and if $1< \lfloor \frac{N_X(T)}{\widetilde{K}T + 1} \rfloor $ we apply the inductive step to obtain a parametrized periodic orbit $(\gamma_2,T_2\leq T)$ such that $S^T_2 = \{(\gamma_1,T_1),(\gamma_2,T_2\leq T) \}$ satisfies (a) and (b). We can go on applying the inductive step to produce sets $S^T_m = \{(\gamma_i,T_i); 1\leq i \leq m\}$ satisfying the desired conditions (a) and (b) as long as $m-1$ is smaller than $\lfloor \frac{N_X(T)}{\widetilde{K}T +1} \rfloor $. By this process we can construct a set $S^T_X =\{(\gamma_i,T_i); 1 \leq i \leq \lfloor \frac{N_X(T)}{\widetilde{K}T +1 } \rfloor \}$ such that for all $i,j \in \{1,..., \lfloor \frac{N_X(T)}{\widetilde{K}T + 1} \rfloor\}$ (a) and (b) above hold true.

For each $i \in \{1,..., \lfloor \frac{N_X(T)}{\widetilde{K}T + 1} \rfloor\}$ let $q_i = \gamma_i(0)$. We define the set $P^T_X:=\{q_i ; 1 \leq i \leq \lfloor \frac{N_X(T)}{\widetilde{K}T } \rfloor  + 1 \}$. The condition (b) satisfied by $S^T_X $ implies that $P^T_X$ is a $T,\delta$-separated set.
It then follows from the definition of the $\delta$- entropy $h_\delta$ that

\begin{equation}
h_{\delta}(\phi_X)\geq \limsup_{T \to +\infty } \frac{\log (\lfloor \frac{N_X(T)}{\widetilde{K}T + 1 } \rfloor)}{T}.
\end{equation}

\

\textbf{Step 5:} Suppose now that for constants $a>0$ and $b$ we have $N_X(T) \geq e^{aT+b}$.

For every $\epsilon >0$ we know that for $T$ big enough we have $e^{\epsilon T}> \widetilde{K}T + 1 $. This implies that
\begin{equation}
\limsup_{T \to +\infty } \frac{\log (\lfloor \frac{N_X(T)}{\widetilde{K}T +1 } \rfloor)}{T} \geq \limsup_{T \to +\infty } \frac{\log (\lfloor \frac{e^{aT+b}}{e^{\epsilon T}} \rfloor)}{T} = \limsup_{T \to +\infty } \frac{\log (\lfloor e^{(a-\epsilon)T+b} \rfloor)}{T}.
\end{equation}

It is clear that $\limsup_{T \to +\infty } \frac{\log (\lfloor e^{(a-\epsilon)T+b} \rfloor)}{T}= a-\epsilon$. We have thus proven that if for constants $a>0$ and $b$ we have $N_X(T) \geq e^{aT+b}$ then $h_{\delta}(\phi_X) \geq a-\epsilon$. Because $\epsilon$ can be taken arbitrarily small we obtain:
\begin{equation}
h_{\delta}(\phi_X) \geq a.
\end{equation}

\

\textbf{Step 6:} We have so far concluded that for all $\delta < \frac{\epsilon_g}{10^6}$ we have $h_{\delta}(\phi_X)\geq a$. We then have:
\begin{equation}
h_{top}(\phi_X) = \lim_{\delta \to 0} h_{\delta}(\phi_X) \geq a,
\end{equation}
finishing the proof of the theorem. \qed

\

\textbf{Remark:} One could naively believe that there exists a constant $\delta_g>0$ depending only on the metric $g$, such that if two parametrized closed curves $\sigma_1:\mathbb{R} \to M$ of period $T_1$ and $\sigma_2:\mathbb{R} \to M$ of period $T_2$ satisfy $\sup_{t\in[0, \max\{T_1,T_2\}] } \{d_g(\sigma_1(t), \sigma_2(t))\} < \delta_g$ then $(\gamma_1,T_1)$ and $(\gamma_2,T_2)$ are freely homotopic to each other. This would make the proof of Theorem \ref{theorem1'} much shorter. However such a constant does not exist. One can easily find for any $\delta >0$ two parametrized curves in the 3-torus which are in different primitive free homotopy classes and satisfy $\sup_{t\in[0, \max\{T_1,T_2\}] } \{d_g(\sigma_1(t), \sigma_2(t))\} < \delta$. We sketch the construction below.

\begin{center}
\begin{minipage}[b]{0.55\textwidth}
    \includegraphics[width=\textwidth]{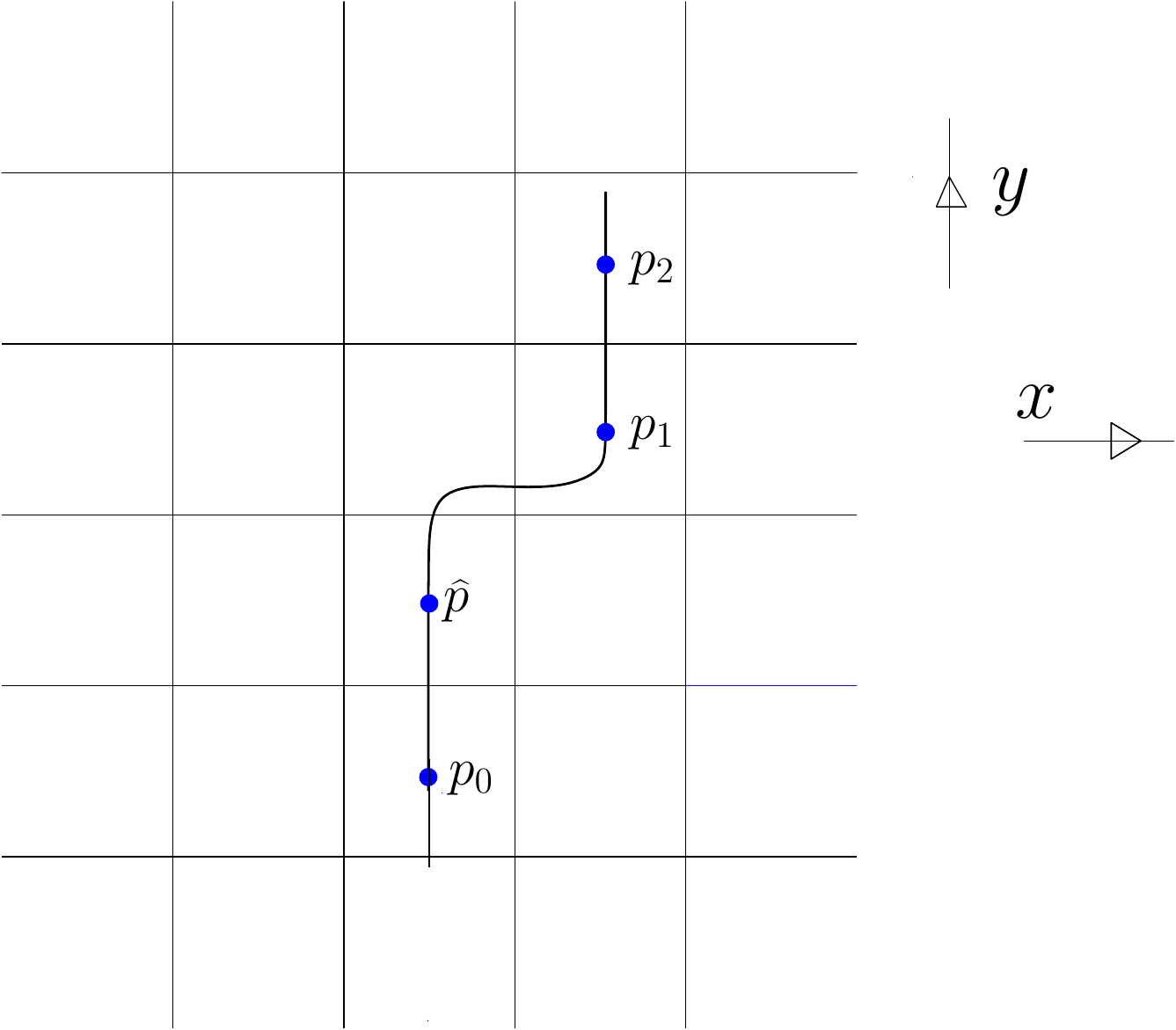}
    \label{fig:1}
    \begin{center}
    Figure 2
    \end{center}
\end{minipage}
\end{center}

Consider coordinates $(x,y,z) \in (\mathbb{R}/\mathbb{Z})^3$ on the three dimensional torus $\mathbb{T}^3$.  Figure 2 above represents the universal cover of the two dimensional torus $\mathbb{T}^2  \subset \mathbb{T}^3 $ obtained by fixing the coordinate $z=0$ in $\mathbb{T}^3$. The dotted points $p_0$, $\widehat{p}$, $p_1$ and $p_2$ in the figure represent lifts of a  point $p \in \mathbb{T}^2$. It is then clear that the curve $c$ represented in the figure projects to a smooth immersed curve in $\mathbb{T}^2 \subset \mathbb{T}^3$.

We consider a parametrization by arc length $\varsigma_1 :[0,T_1] \to \mathbb{R}^2$ of the piece of $c$ connecting $p_0$ and $p_1$. We can extend $\varsigma_1$ periodically to $\mathbb{R}$ by demanding that  $\varsigma_1 (t) = \varsigma_1 (t) + (1,2)$ for all $t \in \mathbb{R}$. This extension is a lift to $\mathbb{R}^2$ of the closed immersed curve obtained by projecting $\varsigma_1([0,T_1])$ to $\mathbb{T}^2$.
 By a very small perturbation of the projection of  $\varsigma_1([0,T_1])$ we can produce a closed smooth embedded curve $\sigma_1 :[0,T_1] \to \mathbb{T}^3$ which closes at the point $(p,0)= \sigma_1(0) = \sigma_1(T_1)$. We consider the natural extension of $\sigma_1 $ to $\mathbb{R}$ obtained by demanding that $\sigma_1(t)= \sigma_1(t - T_1 )$ for all $t \in \mathbb{R}$.

Analogously we consider a parametrization by arc length $\varsigma_2 :[0,T_1 + 1] \to \mathbb{R}^2$ of the piece of $c$ connecting $p_0$ and $p_2$. We can also extend $\varsigma_2$ periodically to $\mathbb{R}$, this time  demanding that  $\varsigma_2 (t) = \varsigma_2  (t) + (1,3)$.
By making a very small perturbation of  $\varsigma_2$ we can produce a closed smooth embedded curve $\sigma_2 :[0,T_1 + 1] \to \mathbb{T}^3$ which closes at the point $(p,\frac{\delta}{K}) =  \sigma_2 (0) = \sigma_2(T_1 + 1)$ and which is disjoint from the image of $\sigma_1$ .We consider the natural extension of $\sigma_2 $ to $\mathbb{R}$ obtained by demanding that $\sigma_2(t)= \sigma_2(t - (T_1+1) )$ for all $t \in \mathbb{R}$.

We point out  that the extensions $\varsigma_1: \mathbb{R} \to \mathbb{R}^2 $ and $\varsigma_2: \mathbb{R} \to \mathbb{R}^2 $ coincide on the interval $[0,T_1 +1]$. To see this just notice that the piece of $c$ connecting $p_0$ and $\widehat{p}$ and the piece of $c$ connecting $p_1$ and $p_2$ project to the same circle in $\mathbb{T}^2$.

Let now $\sigma_0 : [0,T_1 +1] \to \mathbb{T} ^ 2$ be the parametrized curve obtained by projecting $ \varsigma_1 : [0,T_1 +1] \to  \mathbb{R}^2$, which equals $ \varsigma_2 : [0,T_1 +1] \to  \mathbb{R}^2$, to the torus $\mathbb{T}^2$. The curves ${\sigma_{1}}_{|_{[0,T_1 + 1]}}$ and ${\sigma_{2}}_{|_{[0,T_1 + 1]}}$ are both perturbations of the parametrized curve  $\sigma_0$. By making the perturbations sufficiently small we can guarantee that ${\sigma_{1}}_{|_{[0,T_1 + 1]}}$ and ${\sigma_{2}}_{|_{[0,T_1 + 1]}}$ are arbitrarily close. It is immediate to see that ${\sigma_{1}}_{|_{[0,T_1 + 1]}}$ and ${\sigma_{2}}_{|_{[0,T_1 + 1]}}$ are in distinct homotopy classes.

\section{Contact homology} \label{section3}

\subsection{Pseudo-holomorphic curves in symplectic cobordisms}

To define the contact homology theories used in this paper we use pseudoholomorphic curves in symplectizations of contact manifolds and symplectic cobordisms. Pseudo-holomorphic curves were introduced in symplectic manifolds by Gromov in \cite{Gr} and adapted to symplectizations and symplectic cobordisms by Hofer \cite{H}; see also \cite{CPT} as a general reference for pseudoholomorphic curves in symplectic cobordisms.

\subsubsection{Cylindrical almost complex structures}

Let $(Y,\xi)$ be a contact manifold and $\lambda$ a contact form on $(Y,\xi)$. The symplectization of $(Y,\xi)$ is the product $\mathbb{R} \times Y$ with the symplectic form $d(e^s \lambda)$ (where $s$ denotes the $\mathbb{R}$ coordinate in $\mathbb{R} \times Y$). $d\lambda$ restricts to a symplectic form on the vector bundle $\xi$ and it is well known that the set  $\mathfrak{j}(\lambda)$ of $d\lambda$-compatible almost complex structures on the symplectic vector bundle $\xi$ is non-empty and contractible. Notice, that if $Y$ is 3-dimensional the set $\mathfrak{j}(\lambda)$ doesn't depend on the contact form $\lambda$ on $(Y,\xi)$.

For $j \in \mathfrak{j}(\lambda)$ we can define an $\mathbb{R}$-invariant almost complex structure $J$ on $\mathbb{R} \times Y$ by demanding that:

\begin{equation}
J \partial_s = X_\lambda, \ \ J\mid_\xi = j
\end{equation}
We will denote by $\mathcal{J}(\lambda)$ the set of almost complex structures in $\mathbb{R} \times Y$ that are $\mathbb{R}$-invariant, $d(e^s\lambda)$-compatible and satisfy the equation (16) for some $j \in \mathfrak{j}(\lambda)$.

\subsubsection{Exact symplectic cobordisms with cylindrical ends}

An exact symplectic cobordism is, intuitively, an exact symplectic manifold $(W,\varpi)$ that outside a compact subset is like the union of  cylindrical ends of symplectizations. We restrict our attention to exact symplectic cobordisms having only one positive end and one negative end.

Let $(W ,\varpi = d\kappa )$  be an exact symplectic manifold without boundary, and let $(Y^+,\xi^+)$ and $(Y^-,\xi^-)$ be contact manifolds with contact forms $\lambda^+$ and $\lambda^-$. We say that $(W ,\varpi = d\kappa )$ is an exact symplectic cobordism from  $\lambda^+$ to $\lambda^-$ when there exist subsets $W^-$, $W^+$ and $\widehat{W}$ of $W$ and diffeomorphisms $\Psi^+: W^+ \to [0,+\infty) \times Y^+$ and $\Psi^-: W^- \to (-\infty,0] \times Y^-$, such that:
\begin{eqnarray}
\widehat{W} \mbox{ is compact, } W= W^+ \cup \widehat{W} \cup W^- \mbox{ and } W^+ \cap W^- = \emptyset,
\end{eqnarray}
\begin{equation*}
(\Psi^+)^* (e^s\lambda^+) = \kappa \mbox{ and } (\Psi^-)^* (e^s\lambda^-) = \kappa
\end{equation*}

In such a cobordism, we say that an almost complex structure $\overline{J}$ is cylindrical if:
\begin{eqnarray}
\overline{J} \ \mbox{coincides with } \ J^+ \in \mathcal{J}(C^+ \lambda^+) \ \mbox{in the region} \ W^+ \\
\overline{J} \ \mbox{coincides with } \ J^- \in \mathcal{J}(C^- \lambda^-) \ \mbox{in the region} \ W^- \\
\overline{J} \ \mbox{is compatible with} \ \varpi \ \mbox{in} \ \widehat{W}
\end{eqnarray}
where $C^+>0$ and $C^->0$ are constants.

For fixed $J^+ \in \mathcal{J}(C^+ \lambda^+)$ and $J^- \in \mathcal{J}(C^- \lambda^-)$, we denote by $\mathcal{J}(J^-,J^+)$
 set of cylindrical almost complex structures in $(\mathbb{R} \times Y,\varpi)$ coinciding with $J^+$ on $W^+$ and $J^-$ on $W^-$.  It is well known that $\mathcal{J}(J^-,J^+)$ is non-empty and contractible. We will write $ \lambda^+ \succ_{ex} \lambda^-$ when there exists an exact symplectic cobordism from $ \lambda^+$ to $ \lambda^-$ as above. We remind the reader that $ \lambda^+ \succ_{ex} \lambda$ and $ \lambda \succ_{ex} \lambda^-$ implies $ \lambda^+ \succ_{ex} \lambda^-$; or in other words that the exact symplectic cobordism relation is transitive; see \cite{CPT} for a detailed discussion on symplectic cobordisms with cylindrical ends. Notice that a symplectization is a particular case of an exact symplectic cobordism.

\

\textbf{Remark:} we point out to the reader that in many references in the literature, a slightly different definition of cylindrical almost complex structures is used: instead of demanding that $\overline{J}$ satisfies equations (18) and (19), the stronger condition that $\overline{J} \ \mbox{coincides with } \ J^{\pm} \in \mathcal{J}( \lambda^{\pm}) \ \mbox{in the region} \ W^{\pm}$ is demanded. We need to consider this more relaxed definition of cylindrical almost complex structures when we study the cobordism maps of cylindrical contact homologies in subsection~\ref{sec3.3}.

\subsubsection{Splitting symplectic cobordisms}

Let $\lambda^+$, $\lambda$ and $\lambda^-$ be contact forms on $(Y,\xi)$ such that $\lambda^+ \succ_{ex} \lambda$, $ \lambda \succ_{ex} \lambda^-$. For $\epsilon > 0$ sufficiently small, it is easy to see that one also has $\lambda^+ \succ_{ex} (1 + \epsilon)\lambda$ and $ (1-\epsilon)\lambda \succ_{ex} \lambda^-$. Then, for each $R>0$, it is possible to construct an exact symplectic form $\varpi_R = d\kappa_R$ on $W = \mathbb{R} \times Y$ where:
\begin{eqnarray}
\kappa_R = e^{s-R-2}\lambda^+ \ \mbox{in} \ [R+2, + \infty) \times Y, \\
\kappa_R = f(s)\lambda \ \mbox{in} \ [-R,R] \times Y, \\
\kappa_R = e^{s+R+2}\lambda^- \ \mbox{in} \ (-\infty, -R-2] \times Y,
\end{eqnarray}
and $f: [-R,R]  \to [1-\epsilon,1+\epsilon]$, satisfies $f(-R) = 1-\epsilon$, $f(R) = 1+\epsilon$ and $f'>0$.
In $ (\mathbb{R} \times Y, \varpi_R)$ we consider a compatible cylindrical almost complex structure $\widetilde{J}_R$; but we demand an extra condition on $\widetilde{J}_R$:
\begin{eqnarray}
\widetilde{J}_R \ \mbox{coincides with} \ J \in \mathcal{J}(\lambda) \ \mbox{in} \ [-R,R] \times Y.
\end{eqnarray}

Again we divide $W$ in regions:
$W^+ = [R+2, + \infty) \times Y$, $W(\lambda^+,\lambda)= [ R, R+2] \times Y$, $W(\lambda)=[-R, R] \times Y$, $W(\lambda,\lambda^-)=[-R-2,-R] \times Y$ and $W^-=(-\infty, -R-2] \times Y$. The family of exact symplectic cobordisms with cylindrical almost complex structures $(\mathbb{R} \times Y, \varpi_R, \widetilde{J}_R)$ is called a splitting family from $\lambda^+$ to $\lambda^-$ along $\lambda$.

\subsubsection{Pseudoholomorphic curves}

Let $(S,i)$ be a closed Riemann surface without boundary, $\Gamma \subset S$ be a finite set.
Let $\lambda$ be a contact form in $(Y,\xi)$ and $J \in \mathcal{J}(\lambda)$. A finite energy pseudoholomorphic curve in the symplectization $(\mathbb{R} \times Y,J)$ is a map $\widetilde{w}= (r,w):S \setminus \Gamma \to \mathbb{R} \times Y$ that satisfies

\begin{eqnarray}
\overline{\partial}_J(\widetilde{w})= d\widetilde{w} \circ i - J \circ d\widetilde{w}=0
\end{eqnarray}
and
\begin{equation}
0<E(\widetilde{w})= \sup_{q \in \mathcal{E}} \int_{S \setminus \Gamma} \widetilde{w}^*d(q\lambda)
\end{equation}
where $\mathcal{E}= \{ q: \mathbb{R} \to [0,1]; q' \geq 0\}$. The quantity $E(\widetilde{w})$ is called the Hofer energy and was introduced in \cite{H}. The operator $\overline{\partial}_J$ above is called the Cauchy-Riemann operator for the almost complex structure $J$.

For an exact symplectic cobordism $(W,\varpi)$ from $\lambda^+$ to $\lambda^-$ as considered above, and $\overline{J} \in \mathcal{J}(J^-,J^+)$ a finite energy pseudoholomorphic curve is again a map $\widetilde{w}: (S \setminus \Gamma \to W$ satisfying:

\begin{eqnarray}
d\widetilde{w} \circ i = \overline{J} \circ d\widetilde{w},
\end{eqnarray}
and
\begin{equation}
0<E_{\lambda^-}(\widetilde{w}) + E_c(\widetilde{w}) + E_{\lambda^+} (\widetilde{w}) < +\infty,
\end{equation}
where:
\\
$E_{\lambda^-}(\widetilde{w}) = \sup_{q \in \mathcal{E}} \int_{\widetilde{w}^{-1}(W^-))} \widetilde{w}^*d(q\lambda^-)$,
\\
$E_{\lambda^+}(\widetilde{w}) = \sup_{q \in \mathcal{E}} \int_{\widetilde{w}^{-1}(W^+)} \widetilde{w}^*d(q\lambda^+)$,
\\
$E_c(\widetilde{w}) = \int_{\widetilde{w}^{-1}W(\lambda^-,\lambda^+)} \widetilde{w}^*\varpi$.
\\
These energies were also introduced in \cite{H}.

In splitting symplectic cobordisms we use a slightly modified version of energy. Instead of demanding $0<E_-(\widetilde{w}) + E_c(\widetilde{w}) + E_+ (\widetilde{w}) < +\infty$ we demand:

\begin{equation}
0<E_{\lambda^-}(\widetilde{w}) + E_{\lambda^-,\lambda}(\widetilde{w}) + E_{\lambda}(\widetilde{w})+ E_{\lambda,\lambda^+}(\widetilde{w})+ E_{\lambda^+} (\widetilde{w}) < +\infty
\end{equation}
where:
\\
$ E_{\lambda}(\widetilde{w})= \sup_{q \in \mathcal{E}} \int_{\widetilde{w}^{-1}W(\lambda)} \widetilde{w}^*d(q\lambda)$,
\\
$E_{\lambda^-,\lambda}(\widetilde{w})= \int_{\widetilde{w}^{-1}(W(\lambda^-,\lambda))} \widetilde{w}^*\varpi$, \
\\
$E_{\lambda,\lambda^+}(\widetilde{w}=\int_{\widetilde{w}^{-1}(W(\lambda,\lambda^+))} \widetilde{w}^*\varpi$,
\\
and $E_{\lambda^-}(\widetilde{w})$ and $E_{\lambda^+}(\widetilde{w})$ are as above.

\

The elements of the set $\Gamma \subset S$ are called punctures of the pseudoholomorphic $\widetilde{w}$. The work of Hofer et al. \cite{H, P1} allows us do classify the punctures in two types: positive punctures and negative punctures. This classification is done according to the behaviour of $\widetilde{w}$ in the neighbourhood of the puncture.
Before presenting this classification we introduce some notation. Let $B_\delta (z)$ be the ball of radius $\delta$ centered at the puncture $z$, and denote by $\partial (B_\delta (z))$ its boundary. With this in hand, we can describe the types of punctures as follows:
\begin{itemize}
\item{$z \in \Gamma$ is called positive interior puncture when $z \in \Gamma$ and $\lim_{z' \to z} s(z') = +\infty$, and there exist a sequence $\delta_n \to 0$ and  Reeb orbit $\gamma^+$ of  $X_{\lambda^+}$, such that $w(\partial (B_{\delta_n} (z)))$ converges in $C^\infty $ to $\gamma^+$ as $n\to +\infty$}
\item{$z \in \Gamma$ is called negative interior puncture when $z \in \Gamma$ and $\lim_{z' \to z} s(z') = -\infty$, and there exist a sequence $\delta_n \to 0$ and  Reeb orbit $\gamma^-$ of  $X_{\lambda^-}$, such that $w(\partial (B_{\delta_n} (z)))$ converges in $C^\infty $ to $\gamma^-$ as $n\to +\infty$.}
\end{itemize}
The results in \cite{H} and \cite{P1}  imply that these are indeed the only real possibilities we need to consider for the behaviour of the $\widetilde{w}$ near punctures. Intuitively, we have that at the punctures, the pseudoholomorphic curve $\widetilde{w}$ detects Reeb orbits. When for a puncture $z$, there is a subsequence $\delta_n$ such that  $w(\partial(B_{\delta_n} (z)))$ converges to a Reeb orbit $\gamma$, we will say that $\widetilde{w}$ is asymptotic to this Reeb orbit $\gamma$ at the puncture $z$.

If a pseudoholomorphic curve is asymptotic to a non-degenerate Reeb orbit at a puncture, more can be said about its asymptotic behaviour in neighbourhoods of this puncture.
In order to describe the behaviour of $\widetilde{w}$ near a puncture $z$, we take a neighbourhood $U \subset S$ of $z$ that admits a holomorphic chart $\psi_U : (U,z) \to (\mathbb{D},0)$. Using polar coordinates $(r,t) \in (0,+\infty) \times S^1$ we can write $x \in (\mathbb{D} \setminus 0)$ as $x = e^{-r}t $.
With this notation, it is shown in \cite{H} \cite{P1}, that if $z$ is a positive interior puncture on which $\widetilde{w}$ is asymptotic to a non-degenerate Reeb orbit $\gamma^+$ of $X_{\lambda^+}$, then $\widetilde{w} \circ \psi_u^{-1} (r,t) = (s(r,t),w(r,t))$ satisfies:
\begin{itemize}
\item{$w^r(t)=w(r,t)$ converges uniformly in $C^{\infty}$ to a Reeb orbit $\gamma^+$ of $X_{\lambda^+}$, and the convergence rate is exponential.}
\end{itemize}
Similarly, if $z$ is a negative interior puncture on which $\widetilde{w}$ is asymptotic to a non-degenerate Reeb orbit $\gamma^-$ of $X_{\lambda^-}$, then $\widetilde{w} \circ \psi_u^{-1} (r,t) = (s(r,t),w(r,t))$ satisfies:
\begin{itemize}
\item{$w^r(t)=w(r,t)$ converges uniformly in $C^{\infty}$ to a Reeb orbit $\gamma^-$ of $-X_{\lambda^-}$ as  $r \to +\infty$, and the convergence rate is exponential.}
\end{itemize}

\textit{Remark: the fact that the convergence of pseudoholomorphic curves near punctures to Reeb orbits is of exponential nature is a consequence of the asymptotic formula obtained in \cite{P1}. Such formulas are necessary for the Fredholm theory that gives the dimension of the space of pseudoholomorphic curves with fixed asymptotic data.}

The discussion above can be summarised by saying that near punctures the finite pseudoholomorphic curves detect Reeb orbits. It is exactly this behavior that makes these objects useful for the study of dynamics of Reeb vector fields.

For us it will be important to consider the moduli spaces $\mathcal{M}(\gamma,\gamma'_1,...,\gamma'_m;J)$ of of genus $0$ pseudoholomorphic curves, modulo biholomorphic reparametrisation,  with one positive puncture asymptotic to a non-degenerate Reeb orbit $\gamma$ and negative punctures asymptotic to non-degenerate orbits $\gamma'_1,...,\gamma'_m$. It is well known that the linearization $D\overline{\partial}_J$ and $\partial_J$ at any element $\mathcal{M}(\gamma,\gamma'_1,...,\gamma'_m;J)$ is a Fredholm map (we remark that this property is valid for more general moduli spaces of curves with prescribed asymptotic behaviour). One would like to conclude that the dimension of a connected component of $\mathcal{M}(\gamma,\gamma'_1,...,\gamma'_m;J)$ is given by the Fredholm index of an element of $\mathcal{M}(\gamma,\gamma'_1,...,\gamma'_m;J)$. However this is not always the case as problems might appear when the moduli space contain multiply covered pseudoholomorphic curves.

\textbf{Fact: as a consequence of the exactness of the symplectic cobordisms considered above we obtain that the energy $E(\widetilde{w})$ of $\widetilde{w}$ satisfies $E(\widetilde{w}) \leq 5A(\widetilde{w})$ where $ A(\widetilde{w})$ is the sum of the action of the Reeb orbits detected by the punctures of $\widetilde{w}$ counted with multiplicity.}

\subsection{Contact homologies}

Contact homologies were introduced in \cite{SFT} as homology theories which are topological invariants of contact manifolds. In subsections \ref{sec3.2.1} and \ref{sec3.2.2} we give an introduction to the more basic and well known versions of contact homologies. This serves mainly as a motivation to section \ref{sec3.3} where we define the version of contact homology that will be used in this paper. \footnote{We stress that while the versions of contact homology presented in subsections \ref{sec3.2.1} and \ref{sec3.2.2} do depend on the Polyfolds technology currently being developed by Hofer, Wysocki and Zehnder, the version of contact homology which we use in this paper and is presented in subsection \ref{sec3.3} \textbf{does not} depend on Polyfolds and can be constructed in complete rigor with technology that is available in the literature. See the detailed discussion in subsection \ref{sec3.3} below.}

\subsubsection{Full contact homology} \label{sec3.2.1}

Full contact homology was introduced in \cite{SFT} as an important invariant of contact structures. We refer the reader to \cite{SFT} and \cite{B} for detailed presentations of the material contained in this subsection.

Let $(Y^{2\mathrm{n}+1},\xi)$ be a contact manifold with $\lambda$ a non-degenerate contact form. We denote by $\mathcal{P}(\lambda)$ the set of good periodic orbits of the Reeb vector field $X_\lambda$. To each orbit $\gamma \in \mathcal{P}(\lambda)$, we define a $\mathbb{Z}_2$-grading $\mid \gamma \mid = (\mu_{CZ}(\gamma) +(\mathrm{n}-2))\mod2$. An orbit $\gamma$ is called good if it is either simple, or if $\gamma = (\gamma')^i$ for a simple orbit $\gamma'$ with the same grading of $\gamma$.

$\mathfrak{A}(Y,\lambda)$ is defined to be the supercomutative, $\mathbb{Z}_2$ graded, $\mathbb{Q}$ algebra with unit generated by $\mathcal{P}(\lambda)$ (an algebra with this properties is sometimes referred in the literature as a commutative super-algebra or a super-ring). The $\mathbb{Z}_2$-grading on the elements of the algebra is obtained by considering on the generators the grading mentioned above and extending it to $\mathfrak{A}(Y,\lambda)$.

$\mathfrak{A}(Y,\lambda)$ can be equipped with a differential $d_J$. Denote by $\mathcal{M}^k(\gamma,\gamma'_1,...,\gamma'_m;J)$ the moduli space of finite energy pseudoholomorphic curves of genus $0$ and Fredholm index $k$ modulo reparametrization, with one positive puncture asymptotic to $\gamma$ and negative punctures asymptotic to $\gamma'_1,...,\gamma'_m$ in the symplectization $(\mathbb{R} \times Y,J)$. As the almost complex structure $J$ is $\mathbb{R}$-invariant in $\mathbb{R} \times Y$, we have an $\mathbb{R}$-action on $\mathcal{M}^k(\gamma,\gamma'_1,...,\gamma'_m;J)$ and we denote by $\widehat{\mathcal{M}}^k(\gamma,\gamma'_1,...,\gamma'_m;J)= \mathcal{M}^k(\gamma,\gamma'_1,...,\gamma'_m;J)/_{\mathbb{R}}$. Lastly we denote by $\overline{\mathcal{M}}^k(\gamma,\gamma'_1,...,\gamma'_m;J)$ the compactification of $\widehat{\mathcal{M}}^k(\gamma,\gamma'_1,...,\gamma'_m;J)$ as presented in \cite{CPT}. The moduli space $\overline{\mathcal{M}}^k(\gamma,\gamma'_1,...,\gamma'_m;J)$ also involves pseudoholomorphic buildings that appear as limits of a sequence of curves in $\widehat{\mathcal{M}}^k(\gamma,\gamma'_1,...,\gamma'_m;J)$  that ``breaks''; we refer the reader to \cite{CPT} for a more detailed description of these moduli spaces.  To define our differential we need the following hypothesis:

\textbf{Hypothesis H}: there exists an abstract perturbation of the Cauchy-Riemann operator $\partial_J $ such that the compactified moduli spaces $\overline{\mathcal{M}}(\gamma,\gamma'_1,...,\gamma'_m;J)$ of solutions of the perturbed equation are unions of branched manifolds with corners and rational weights whose dimension is given by the Conley-Zehnder index of the asymptotic orbits and the relative homology class of the solution.

The proof that Hypothesis H is true is still not written. Establishing its validity is one of the main reasons for the development of the Polyfold technology by Hofer, Wysocki and Zehnder. We define:

\begin{equation}
d_J \gamma = m(\gamma) \sum_{\gamma'_1,...,\gamma'_m} \frac{C(\gamma,\gamma'_1,...,\gamma'_m)}{m!} \gamma'_1 \gamma'_2... \gamma'_m
\end{equation}
where $C(\gamma,\gamma'_1,...,\gamma'_m)$ is the algebraic count of points in the $0$-dimensional manifold
\begin{equation}
\widehat{\mathcal{M}}^1(\gamma,\gamma'_1,...,\gamma'_2;J)
\end{equation}
and $m(\gamma)$ is the multiplicity of $\gamma$. $d_J$ is extended to the whole algebra by the Leibnitz rule. Under hypothesis \textbf{H} it was proved in \cite{SFT} that $(d_J)^2 =0$. We have therefore that $(\mathfrak{A}(Y,\lambda),d_J)$ is a differential $\mathbb{Z}_2$ graded super-commutative algebra. We define:

\begin{definition}
\ The \textbf{full contact homology} $C\mathbb{H}(\lambda,J)$ of $\lambda$ is the homology of the complex $(\mathfrak{A},d_J)$.
\end{definition}

Under Hypotesis H, it was also proved in \cite{SFT} that the full contact homology does not depend on the contact form $\lambda$ on $(Y,\xi)$ nor on the choice of the cylindrical almost complex structure $J\in \mathcal{J}(\lambda)$.

\subsubsection{Cylindrical contact homology} \label{sec3.2.2}
Suppose now that $(Y,\xi)$ is a contact manifold, and $\lambda$ is a non-degenerate hypertight contact form on $(Y,\xi)$. Fix a cylindrical almost complex structure $J \in \mathcal{J}(\lambda)$. For hypertight contact manifolds we can define a simpler version of contact homology called cylindrical contact homology. We denote by $CH_{cyl}(\lambda)$ the $\mathbb{Z}_2$-graded $\mathbb{Q}$-vector space generated by the elements of $\mathcal{P}(\lambda)$. The differential $d^{cyl}_J:CH_{cyl}(\lambda) \to CH_{cyl}(\lambda)$ will count elements in the moduli space $\widehat{\mathcal{M}}^1(\gamma,\gamma';J)$. For the generators $\gamma \in \mathcal{P}(\lambda) $ we define
\begin{equation}
d^{cyl}_J(\gamma)=cov(\gamma)\sum_{\gamma' \in \mathcal{P}(\lambda)} C(\gamma,\gamma';J)\gamma',
\end{equation}
where $C(\gamma,\gamma';J)$ is the algebraic count of elements in $\widehat{\mathcal{M}}^1(\gamma,\gamma';J)$ and  $cov(\gamma)$ is the covering number of $\gamma$. For $\lambda$ hypertight and assuming Hypothesis H is true, Eliashberg, Givental and Hofer proved in \cite{SFT} that $(d^{cyl}_J)^2=0$.

\begin{definition}
The \textbf{cylindrical contact homology} $C\mathbb{H}_{cyl}(\lambda)$ of $\lambda$ is the homology of the complex $(CH_{cyl}(\lambda),d^{cyl}_J)$.
\end{definition}
Under Hypothesis H, the cylindrical contact homology doesn't depend on the hypertight contact form $\lambda$ on $(Y,\xi)$ nor on the cylindrical almost complex structure $J\in \mathcal{J}(\lambda)$.

Denote by $\Lambda$ the set of free homotopy classes of $Y$. It is easy to see that for each $\rho \in \Lambda$ the subspace $CH_{cyl}^{\rho}(\lambda) \subset CH_{cyl}(\lambda)$ generated by the set $\mathcal{P}_{\rho}(\lambda)$ of good periodic orbits in $\rho$ is a subcomplex of $(CH_{cyl}(\lambda),d^{cyl}_J)$. This follows from the fact that the number $C(\gamma,\gamma';J)$ can only be non-zero for Reeb orbits $\gamma'$ that are freely homotopic to $\gamma$, which implies that the restriction $d^{cyl}_J\mid_{CH_{cyl}^{\rho}(\lambda)}$ has image in $CH_{cyl}^{\rho}(\lambda)$. From now on we will denote the restriction $d^{cyl}_J\mid_{CH_{cyl}^{\rho}}:CH_{cyl}^{\rho}(\lambda) \to CH_{cyl}^{\rho}(\lambda)$ by $d^{\rho}_J$. Denoting by $C\mathbb{H}_{cyl}^{\rho}$ the homology of $(CH_{cyl}^{\rho}(\lambda),d^{\rho}_{J})$ we thus have
\begin{equation}
C\mathbb{H}_{cyl}(\lambda)= \bigoplus_{\rho \in \Lambda} C\mathbb{H}_{cyl}^{\rho}.
\end{equation}

The fact that we can define partial versions of cylindrical contact homology restricted to certain free homotopy classes will be of crucial importance for us. It will allow us to obtain our results without resorting to Hypothesis H. This is explained in the next subsection.

\subsubsection{Cylindrical contact homology in special homotopy classes}\label{sec3.3}
Maintaining the notation of the previous sections we denote by $(Y,\xi)$ a contact manifold endowed with a hypertight contact form $\lambda$.

Let $\Lambda_0$ denote the set of primitive free homotopy classes of $Y$.  Let $\rho \in \Lambda$ be either an element of $\Lambda_0$, or a free homotopy class which contains only simple Reeb orbits of $\lambda$. Assume that all Reeb orbits in $\mathcal{P}_\rho(\lambda)$ are non-degenerate. By the work of Dragnev \cite{Dr}, we know that there exists a generic subset $\mathcal{J}^{\rho}_{reg}(\lambda)$ of $\mathcal{J}(\lambda)$ such that for all $J\in \mathcal{J}_{reg}(\lambda)$ we have:
\begin{itemize}
\item{for all Reeb orbits $\gamma_1,\gamma_2 \in \rho$, the moduli space of pseudoholomorphic cylinders $\mathcal{M}(\gamma_1,\gamma_2;J)$ is transverse, i.e. the linearized Cauchy-Riemann operator $D\partial_J(\widetilde{w})$ is surjective for all $\widetilde{w} \in \mathcal{M}(\gamma_1,\gamma_2;J);$}
\item{for all Reeb orbits $\gamma_1,\gamma_2 \in \rho$, each connected component $\mathcal{L}$ of the moduli space $\mathcal{M}(\gamma_1,\gamma_2;J)$ is a manifold whose dimension is given by the Fredholm index of any element $\widetilde{w} \in \mathcal{L}$.}
\end{itemize}

In this case, for $J \in \mathcal{J}_{reg}(\lambda)$, we define:
\begin{equation}
d^{\rho}_J(\gamma)= cov(\gamma)\sum_{\gamma' \in \mathcal{P}_{\rho}(\lambda)} C^{\rho}(\gamma,\gamma';J)\gamma' = \sum_{\gamma' \in \mathcal{P}_{\rho}(\lambda)} C^{\rho}(\gamma,\gamma';J)\gamma'
\end{equation}
where $C^{\rho}(\gamma,\gamma';J)$ is the number of points of the moduli space $\widehat{\mathcal{M}}^1(\gamma,\gamma';J)$. The second equality follows from the fact that all Reeb orbits in $\rho$ are simple, which implies $cov(\gamma)=1$.

For $\lambda$ and $\rho$ as above and $J \in \mathcal{J}^{\rho}_{reg}(\lambda)$, the differential $d^{\rho}_{J}: CH_{cyl}^{\rho}(\lambda) \to CH_{cyl}^{\rho}(\lambda)$ is well-defined and satisfies $(d^{\rho}_{J})^2=0$. Therefore, in this situation, we can define the cylindrical contact homology $C\mathbb{H}_{cyl}^{\rho,J}(\lambda)$ without the need of Hypothesis H. Once the transversality for $J$ has been achieved, and using coherent orientations constructed in \cite{BM}, the proof that $d^{\rho}_{J}$ is well-defined and that $(d^{\rho}_{J})^2=0$ is a combination of compactness and gluing, similar to the proof of the analogous result for Floer homology. For the convenience of the reader we sketch these arguments below:

\

\textbf{For $\rho$ as above, $d^{\rho}_J: CH_{cyl}^{\rho}(\lambda) \to CH_{cyl}^{\rho}(\lambda)$ is well-defined, and for every $\gamma \in \mathcal{P}_\rho(\lambda)$ the differential $d^{\rho}_J(\gamma)$  is a finite sum.} \\
The moduli space $\widehat{\mathcal{M}}^1(\gamma,\gamma';J)$ can be non-empty only if $A(\gamma') \leq A(\gamma)$. It then follows from the non-degeneracy of $\lambda$ that for a fixed $\gamma$ the numbers $C^{cyl}(\gamma,\gamma';J)$ can be nonzero for only finitely many $\gamma'$. To see that $C^{cyl}(\gamma,\gamma';J^)$ is finite for every $\gamma' \in \rho$ suppose by contradiction that there is a sequence $\widetilde{w}_i$ of distinct elements of $\widehat{\mathcal{M}}^1 (\gamma,\gamma';J)$. By the SFT compactness theorem \cite{CPT} such a sequence has a convergent subsequence that converges to a pseudoholomorphic building $\widetilde{w}$ which has Fredholm index 1. Because of the hypertightness of $\lambda$, no bubbling can occur and all the levels $\widetilde{w}^1,...,\widetilde{w}^k$ of the building $\widetilde{w}$ are pseudoholomorphic cylinders. As all Reeb orbits of $\lambda$ in $\rho$ are simple, it follows that all these cylinders are somewhere injective pseudoholomorphic curves, and the regularity of $J$ implies that they must all have Fredholm index $\geq 1$. As a result we have $1=I_F (\widetilde{w})= \sum (I_F(\widetilde{w}^l)) \geq k$, which implies $k=1$. Thus $ \widetilde{w} \in \widehat{\mathcal{M}}^1 (\gamma,\gamma';J)$ and is the limit of a sequence of distinct elements of  $\widehat{\mathcal{M}}^1 (\gamma,\gamma';J)$. This is absurd  because $\widehat{\mathcal{M}}^1 (\gamma,\gamma';J)$ is a $0$-dimensional manifold. We thus conclude that the numbers $C^{cyl}(\gamma,\gamma';J^)$ are all finite. \qed

\

\textbf{For $\rho$ as above, $(d^{\rho}_{J})^2=0$.}
If we write
\begin{equation}
d^{\rho}_J \circ d^{\rho}_J(\gamma) = \sum_{\gamma'' \in \mathcal{P}_{\rho}(\lambda)} m_{\gamma,\gamma''}\gamma'',
\end{equation}
we know that $m_{\gamma,\gamma'}$ is the number of two-level pseudo holomorphic buildings $\widetilde{w} =(\widetilde{w^1},\widetilde{w^2})$ such that $\widetilde{w^1} \in \widehat{\mathcal{M}}^1 (\gamma,\gamma';J)$ and $\widetilde{w^2} \in \widehat{\mathcal{M}}^1 (\gamma',\gamma'';J)$, for some $\gamma' \in \mathcal{P}_{\rho}(\lambda)$. Because of transversality of $\widetilde{w^1}$ and $\widetilde{w^2}$ we can perform gluing. This implies that $\widetilde{w}$ is in the boundary of the moduli space $\overline{\mathcal{M}}^2 (\gamma,\gamma'';J)$. Taking a sequence $\widetilde{w}_i$ of elements in $\widehat{\mathcal{M}}^2 (\gamma,\gamma'';J)$ converging to the boundary of $\overline{\mathcal{M}}^2 (\gamma,\gamma'';J)$ and arguing similarly as above, we have that this sequence converges to a pseudoholomorphic building $\widetilde{w}_{\infty}$, whose levels are somewhere injective pseudoholomorphic cylinders. Using that $I_F(\widetilde{w}_{\infty})=2$ we  obtain that $\widetilde{w}_{\infty}$ can have at most 2 levels. As $\widetilde{w}_{\infty}$ is in the boundary of $\overline{\mathcal{M}}^2 (\gamma,\gamma'';J)$ it cannot have only one level, and is therefore a two-level pseudo holomorphic building whose levels have Fredholm index 1. Summing up, $\widetilde{w}_{\infty} =(\widetilde{w}_{\infty}^1,\widetilde{w}_{\infty}^2)$, where $\widetilde{w}^1_{\infty} \in \widehat{\mathcal{M}}^1 (\gamma,\gamma';J)$ and $\widetilde{w}^2_{\infty} \in \widehat{\mathcal{M}}^1 (\gamma',\gamma'';J)$, for some $\gamma' \in \mathcal{P}_{\rho}(\lambda)$.

The discussion above implies that $m_{\gamma,\gamma''}$ is the count with signs of boundary components of the compactified moduli space $\overline{\mathcal{M}}^2 (\gamma,\gamma'';J)$ which is homeomorphic to a one-dimensional manifold with boundary. Because the signs of this count are determined by coherent orientations of $\overline{\mathcal{M}}^2 (\gamma,\gamma'';J)$, it follows that $m_{\gamma,\gamma''}=0$.
\qed

The discussion above gives us the following
\begin{proposition*}
Let $(Y,\xi)$ be a contact manifold with a hypertight contact form $\lambda$. Let $\rho \in \Lambda$ be either an element of $\Lambda_0$, or a free homotopy class which contains only simple Reeb orbits of $\lambda$. Assume that all Reeb orbits in $\mathcal{P}_\rho(\lambda)$ are non-degenerate and pick $J \in \mathcal{J}^{\rho}_{reg}(\lambda)$. Then, $d^{\rho}_{J}$ is well defined and $(d^{\rho}_{J})^2=0$.
\end{proposition*}

\

Exact symplectic cobordisms induce homology maps for the SFT-invariants. We describe how this is done for the version of cylindrical contact homology considered in this section.
Let $(Y^+,\xi^+)$ and $(Y^-,\xi^-)$ be contact manifolds, with hypertight contact forms $\lambda^+$ and $\lambda^-$. Let $(W,\omega)$ be an exact symplectic cobordism from $\lambda^+$ to $\lambda^-$. Assume that $\rho$ is either a primitive free homotopy class or that all the closed Reeb orbits of both $\lambda^+$ and $\lambda^-$ which belong to $\rho$ are simple. Assume moreover that all Reeb orbits of both $\mathcal{P}_\rho(\lambda^+)$ and $\mathcal{P}_\rho(\lambda^-)$ are non-degenerate. Choose almost complex structures $J^+\in \mathcal{J}^{\rho}_{reg}(\lambda^+)$ and $J^-\in \mathcal{J}^{\rho}_{reg}(\lambda^-)$. From the work of Dragnev \cite{Dr} (see also section 2.3 in \cite{Mo}) we know that there is a generic subset $\mathcal{J}^{\rho}_{reg}(J^-,J^+) \in \mathcal{J}(J^-,J^+)$ such that for $\widehat{J} \in \mathcal{J}^{\rho}_{reg}(J^-,J^+)$, $\gamma^+ \in \mathcal{P}_{\rho}(\lambda^+)$ and $\gamma^- \in \mathcal{P}_{\rho}(\lambda^-)$:

\begin{itemize}
\item{all the curves $\widetilde{w}$ in the moduli spaces $\mathcal{M} (\gamma^+,\gamma^-;\widehat{J})$ are Fredholm regular,}
\item{the connected components $\mathcal{V}$ of $\mathcal{M} (\gamma^+,\gamma^-;\widehat{J})$ have dimension equal to the Fredholm index of any pseudoholomorphic curve in $\mathcal{V}$.}
\end{itemize}

In this case we can define a map $\Phi^{\widehat{J}}: CH_{cyl}^{\rho}(\lambda^+) \to CH_{cyl}^{\rho}(\lambda^-)$ given on elements of $\mathcal{P}_{\rho}(\lambda^+)$ by
\begin{equation}
\Phi^{\widehat{J}}(\gamma^+)= \sum_{\gamma^- \in \mathcal{P}_{\rho}(\lambda^-)} n_{\gamma^+,\gamma^-}\gamma^-,
\end{equation}
where $n_{\gamma^+,\gamma^-}$ is the number pseudoholomorphic cylinders with Fredholm index $0$, positively asymptotic to $\gamma^+$ and negatively asymptotic to $\gamma^-$. Using a combination of compactness and gluing (see \cite{B}) one proves that $\Phi^{\widehat{J}} \circ d^{\rho}_{J^+} = d^{\rho}_{J^-} \circ \Phi^{\widehat{J}} $. As a result we obtain a map $\Phi^{\widehat{J}}:C\mathbb{H}_{cyl}^{\rho,J^+}(\lambda^+) \to C\mathbb{H}_{cyl}^{\rho,J^-}(\lambda^-) $ on the homology level.

We study the cobordism map in the following situation: take $(V=\mathbb{R}\times Y,\varpi)$ to be an exact symplectic cobordism from $C\lambda$ to $c\lambda$ where $C>c>0$, and $\lambda$ is a hypertight contact form. Suppose that one can make an isotopy of exact symplectic cobordisms $(\mathbb{R}\times Y, \varpi_t)$ from $C\lambda$ to $c\lambda$, with $\varpi_t$ satisfying $\varpi_0 = \varpi$ and $\varpi_1= d(e^s\lambda_0)$. We consider the space $\widetilde{\mathcal{J}}(J,J)$ of smooth homotopies
\begin{equation}
t \in [0,1]; J_t \in \mathcal{J}(J,J)
\end{equation}
such that $J_0=J_V$, $J_1 \in \mathcal{J}_{reg}(\lambda)$, and $J_t$ is compatible with $\varpi_t$ for every $t\in [0,1]$. $J_t$ is a deformation of $J_0$ to $J_1$, through asymptotically cylindrical almost complex structures in the cobordisms $(\mathbb{R}\times Y,\varpi_t)$. For Reeb orbits $\gamma,\gamma' \in \mathcal{P}_{\rho}(\lambda)$ we consider the moduli space
\begin{equation}
\widetilde{\mathcal{M}}^1(\gamma,\gamma';J_t)= \{(t,\widetilde{w}) \mid t \in [0,1] \mbox{ and } \widetilde{w}\in  \widehat{\mathcal{M}}^1(\gamma,\gamma';J_t) \}.
\end{equation}

By using the techniques of \cite{Dr}, we know that there is a generic subset $\widetilde{\mathcal{J}}_{reg}(J,J)=\widetilde{\mathcal{J}}(J,J)$  such that $\widetilde{\mathcal{M}}^1(\gamma,\gamma';J_t)$ is a 1-dimensional smooth manifold with boundary. The crucial condition that makes this valid is again the fact that the all the pseudoholomorphic curves that make part of this moduli space are somewhere injective.

We have the following proposition which is a consequence of the combination of work of Eliashberg, Givental and Hofer \cite{SFT} and Dragnev \cite{Dr}.

\begin{proposition} \label{proposition1}
Let $(Y,\xi)$ be a contact manifold with a hypertight contact form $\lambda$. Let $\lambda^+ = C\lambda$ and $\lambda^- = c\lambda$ where $C>c>0$ are constants, and $\rho$ be either a primitive free homotopy class or a free homotopy class in which all Reeb orbits of $\lambda$ are simple. Assume that all Reeb orbits in $\mathcal{P}_\rho(\lambda)$ are non-degenerate.  Choose an almost complex structure $J\in \mathcal{J}^{\rho}_{reg}(\lambda)$, and set $J^+= J^- =J$. Let $(W=\mathbb{R} \times Y,\varpi)$ be an exact symplectic cobordism from $C\lambda$ to $c\lambda$, and choose a regular almost complex structure $\widehat{J} \in \mathcal{J}^{\rho}_{reg}(J^-,J^+)$. Then, if there is an homotopy $(\mathbb{R}\times Y,\varpi_t)$ of exact symplectic cobordisms from $C\lambda$ to $c\lambda$, with $\varpi_0 = \varpi$ and $\varpi_1= d(e^s\lambda)$, it follows that the map $\Phi^{\widehat{J}}:CH_{cyl}^{\rho,J}(\lambda) \to CH_{cyl}^{\rho,J}(\lambda) $ is chain homotopic to the identity.
\end{proposition}

The proof is again a combination of compactness and gluing, and we sketch it below. We refer the reader to \cite{B} and \cite{SFT} for the details.

\textit{Sketch of the proof:} We define initially the following map
\begin{equation}
K: CH^{\rho}_{cyl}(\lambda) \to CH^{\rho}_{cyl}(\lambda)
\end{equation}
that counts finite energy, Fredholm index $-1$ pseudoholomorphic cylinders in the cobordisms $(\mathbb{R} \times Y,\varpi_t)$ for $t \in [0,1]$. Because of the regularity of our homotopy, the moduli space of index $-1$ cylinders whose positive puncture detects a fixed Reeb orbit $\gamma$ is finite, and therefore the map $K$ is well defined.

Notice that for $t=1$ the cobordism map $\Phi^{\widehat{J_1}}$ is the identity, and the pseudoholomorphic curves that define it are just trivial cylinders over Reeb orbits. For $t=0$, $\Phi^{\widehat{J_0}}$ counts index $0$ cylinders in the cobordisms $(\mathbb{R} \times Y,\varpi)$. From the regularity of $J_0$, $J_1$ and the homotopy $J_t$, we have that the pseudoholomorphic cylinders involved in these two maps belong to the $1$-dimensional moduli spaces $\widetilde{\mathcal{M}}^1(\gamma,\gamma';J_t)$.

By using a combination of compactness and gluing we can show that the boundary of the moduli space $\widetilde{\mathcal{M}}^1(\gamma,\gamma';J_t)$ is exactly the set of pseudoholomorphic buildings $\widetilde{w}$ with two levels $\widetilde{w}_{cob}$ and $\widetilde{w}_{symp}$ such that: $\widetilde{w}_{cob}$ is an index $-1$ cylinder in a cobordism $(\mathbb{R} \times Y,\varpi_t)$ and $\widetilde{w}_{symp}$ is an index $1$ pseudoholomorphic cylinder in the symplectization of $\lambda$ above or below $\widetilde{w}_{cob}$. Such two level buildings are exactly the ones counted in the map $ K \circ d^{cyl}_J + d^{cyl}_J\circ K$. As a consequence one has that the difference between the maps $\Phi^{\widehat{J_1}}=Id$ and $\Phi^{\widehat{J}}$ is equal to $ K \circ d^{cyl}_J + d^{cyl}_J\circ K$. This implies that $\Phi^{\widehat{J}}$ is chain homotopic to the identity.
\qed

The result above can be used to show that $C\mathbb{H}_{cyl}^{\rho}(\lambda)$ does not depend on the regular almost complex structure $J$ used to define the differential $d_J$.

\section{Exponential homotopical growth rate of $C\mathbb{H}_{cyl}(\lambda_0)$ and estimates for $h_{top}$} \label{section4}

In this section we define the exponential homotopical growth of contact homology and relate it to the topological entropy of Reeb vector fields. The basic idea is to use non-vanishing of cylindrical contact homology of $(M,\xi)$ in a free homotopy class to obtain existence of Reeb orbits in such an homotopy class for any contact form on $(M,\xi)$; this idea is present in \cite{HMS,Mo}. It is straightforward to see that the period and action of a Reeb orbit are equal and in the sequel we will use the same notation to refer period and action of Reeb orbits.

Let $(M,\xi)$ be a contact manifold and $\lambda_0$ be a hypertight contact form on  $(M,\xi)$. For $T>0$ we define $\widetilde{\bigwedge}_T(\lambda_0)$ to be the set of free homotopy classes of $M$ such that $\rho \in \widetilde{\bigwedge}_T(\lambda_0)$ if, and only if, all Reeb orbits of $X_{\lambda_0}$ in $\rho$ are simply covered, non-degenerate, have action/period smaller than $T$ and $C\mathbb{H}_{cyl}^{\rho}(\lambda_0) \neq 0$. We define $N^{cyl}_T(\lambda_0)$ to be the cardinality $\sharp(\widetilde{\bigwedge}_T(\lambda_0))$.

\

\textbf{Definition:} We say that the cylindrical contact homology $C\mathbb{H}_{cyl}(\lambda_0)$ of $(M,\lambda_0)$ has exponential homotopical growth with exponential weight $a>0$ if there exist $T_0\geq 0$ and $b$, such that for all $T\geq T_0$ $N^{cyl}_T(\lambda_0)=\sharp(\widetilde{\bigwedge}_T(\lambda_0))\geq e^{aT + b}$.

\

The main result of this section is the following:

\begin{theorem} \label{theorem2'}
Let $\lambda_0$ be a hypertight contact form on a contact manifold $(M,\xi)$ and assume that the cylindrical contact homology $C\mathbb{H}_{cyl}(\lambda_0)$ has exponential homotopical growth with exponential weight $a>0$. Then for every $C^k$ ($k\geq2$) contact form $\lambda$  on  $(M,\xi)$ the Reeb flow of $X_\lambda$ has positive topological entropy. More precisely, if $f_\lambda$ is the unique function such that $\lambda = f_\lambda \lambda_0$, then
\begin{equation}
h_{top}(X_{\lambda})\geq \frac{a}{\max f_\lambda}.
 \end{equation}
\end{theorem}

\textit{Proof:}
We write $E=\max f$.

\textbf{Step 1:} \\
 We assume initially that $\lambda$ is non-degenerate and $C^{\infty}$.
For every $\epsilon>0$ is possible to construct an exact symplectic cobordism from $(E+\epsilon) \lambda_0$ to $\lambda$. Analogously, for $e>0$ small enough, it is possible to construct an exact symplectic cobordism from $\lambda$ to $e\lambda_0$.

Using these cobordisms, it is possible to construct a splitting family $(\mathbb{R} \times M, \varpi_R, J_R)$ from $(E+\epsilon) \lambda_0$ to $e\lambda_0$, along $\lambda$, such that for every $R>0$ $(\mathbb{R} \times M, \varpi_R, J_R)$ is homotopical to the symplectization of $\lambda_0$. We fix a regular almost complex structure $J_0 \in \mathcal{J}^{\rho}_{reg}(\lambda_0)$ and $J \in \mathcal{J}(\lambda)$, and demand that $J_R$ coincides with $J_0$ in the positive and negative ends of the cobordism, and with $J$ on $[-R,R] \times M$.

Let $\rho \in \widetilde{\bigwedge}_T(\lambda_0)$. We claim that for every $R$ there exists a finite energy pseudoholomorphic cylinder $\widetilde{w}$ in $(\mathbb{R} \times M,J_R)$ positively asymptotic to a Reeb orbit in $\mathcal{P}_\rho(\lambda_0)$ and negatively asymptotic to an orbit in $\mathcal{P}_\rho(\lambda_0)$.

If this was not true for a certain $R>0$, then because of the absence of pseudoholomorphic cylinders asymptotic to Reeb orbits in $\mathcal{P}_\rho(\lambda_0)$ we would have that $J_R \in \mathcal{J}^{\rho}_{reg}(J_0,J_0)$. Therefore, the map $\Phi^{J_R}:C\mathbb{H}_{cyl}^{\rho}(\lambda_0) \to C\mathbb{H}_{cyl}^{\rho}(\lambda_0)$ induced by $(\mathbb{R} \times M, \varpi_R, J_R)$ is well-defined. But because there are no pseudoholomorphic cylinders asymptotic to Reeb orbits in $\mathcal{P}_\rho(\lambda_0)$, we have that the map $\Phi^{J_R}:C\mathbb{H}_{cyl}^{\rho}(\lambda_0) \to C\mathbb{H}_{cyl}^{\rho}(\lambda_0)$ vanishes. On the other hand, from subsection \ref{sec3.3} we know that $\Phi^{J_R}$ the identity. As $\Phi^{J_R}$ vanishes and is
the identity we conclude that $C\mathbb{H}_{cyl}^{\rho}(\lambda_0)=0$, contradicting that $\rho \in \widetilde{\bigwedge}_T(\lambda_0)$.

\

\textbf{Step 2:} \\
Let $\rho \in \widetilde{\bigwedge}_T(\lambda_0) $, $R_n \to +\infty$ be a strictly increasing sequence and $\widetilde{w}_n : \mathbb{R} \times (S^1 \times \mathbb{R},i) \to (\mathbb{R} \times M, J_{R_n}) $ be a sequence of pseudoholomorphic cylinders with one positive puncture asymptotic to an orbit in $\mathcal{P}_\rho(\lambda_0)$ and one negative puncture asymptotic to an orbit in $\mathcal{P}_\rho(\lambda_0)$. Notice that, because of the properties of $\rho$ the energy of $\widetilde{w}_n $ is uniformly bounded.

Therefore we can apply the SFT compactness theorem to obtain a subsequence of $\widetilde{w}_n $ which converges to a pseudoholomorphic buiding $\widetilde{w}$. Notice that in order to apply the SFT compactness theorem we need to use the non-degeneracy of $\lambda$. Moreover we can give a very precise description of the building.

Let $\widetilde{w}^k$ for $k\in \{1,...,m\}$ be the levels of the pseudoholomorphic building $\widetilde{w}$. Because the topology of our curve doesn't change on the breaking we have the following picture:
\begin{itemize}
\item{the upper level $\widetilde{w}^1$ is composed by one connected pseudoholomorphic curve, which has one positive puncture asymptotic to an orbit $\gamma_0 \in \mathcal{P}_\rho(\lambda_0)$, and several negative punctures. All of the negative punctures detect contractible orbits, except one that detects a Reeb orbit $\gamma_1$ which is also in $\rho$.}
\item{on every other level $\widetilde{w}^k$ there is a special pseudoholomorphic curve which has one positive puncture asymptotic to a Reeb orbit $\gamma_{k-1}$ in $\rho$, and at least one but possibly several negative punctures. Of the negative punctures there is one that is asymptotic to an orbit $\gamma_k$ in $\rho$, while all the others detect contractible Reeb orbits.}
\end{itemize}

Because of the splitting behavior of the cobordisms $(\mathbb{R} \times M, J_{R_n}$ it is clear that there exists a $k_0$, such that the level $\widetilde{w}^k$ is in an exact symplectic cobordism from $(E+\epsilon)\lambda_0$ to $\lambda$. This implies that the special orbit $\gamma_{k_0}$ is a Reeb orbit of $X_\lambda$ in the homotopy class $\rho$.

Notice that $A(\gamma_0) \leq (E+\epsilon)T$. This implies that all the other orbits appearing as punctures of the building $\widetilde{w}$ have action smaller than $(E+\epsilon)T$, and in particular that $\gamma_{k_0}$ has action smaller than $(E+\epsilon)T$.

As we can do the construction above for any $\epsilon>0$ we can obtain a sequence of Reeb orbits $\gamma_{K_0(j)}$ which are all in $\rho$ and such that $A(\gamma_{K_0(j)}) \leq (E+\frac{1}{j})T$. Using Arzela-Ascoli's Theorem one can extract a convergent subsequence of $\gamma_j^\rho$. Its limit $\gamma_\rho$ is clearly a Reeb orbit of $\lambda$ in the free homotopy class $\rho$ and with action smaller or equal to $ET$.

\

\textbf{Step 3:} Estimating $N_{X_{\lambda}}(T)$. \footnote{Recall that as we defined in section \ref{section2}, $N_{X_{\lambda}}(T)$ is the number of distinct free homotopy classes of $M$ that contain periodic orbits of $X_\lambda$ with period $\leq T$.}

From step 2, we know that if $\rho \in \widetilde{\bigwedge}_T(\lambda_0)$ then there is a Reeb orbit $\gamma_\rho$ of the Reeb flow of $X_\lambda$ with $A(\gamma_\rho)\leq ET$. Recalling that the period and the action of a Reeb orbit coincide we obtain that $N_{X_{\lambda}}(T) \geq \sharp(\widetilde{\bigwedge}_{\frac{T}{E}}(\lambda_0))$. Under the hypothesis of the theorem we have
\begin{equation}
N_{X_{\lambda}}(T) \geq e^{\frac{aT}{E} + b}.
\end{equation}
Applying Theorem \ref{theorem1'} we then obtain $h_{top}(X_{\lambda}) \geq \frac{a}{E}$. This proves the theorem in the case $\lambda$ is $C^{\infty}$ and non-degenerate.

\

\textbf{Step 4:} Passing to the case of a general $C^{k\geq2}$ contact form $\lambda$ (the case where $\lambda$ is degenerate is included here).

Let $\lambda_i$ be a sequence of non-degenerate contact forms converging in the $C^k$-topology to a contact form $\lambda$ which is $C^k $ ($k\geq 2$) and possibly degenerate. For every $\epsilon >0$ there is $i_0$ such that for $i>i_0$; there exists an exact symplectic cobordism from $(E+\epsilon)\lambda_0 $ to $\lambda_i$.

Fixing then an homotopy class $ \rho \in \widetilde{\bigwedge}_T(\lambda_0)$ we know, by the previous steps, that there exists a Reeb orbit $\gamma_\rho(n)$ of $\lambda_n$ in the homotopy class $\rho$ with action smaller than $(E+\epsilon)T$. By taking the sequence $\gamma_\rho(n)$ and applying Arzela-Ascoli's theorem we obtain a subsequence which converge to a Reeb orbit $\gamma_{\epsilon,\rho}$ of $X_\lambda$ with$ A(\gamma_{\epsilon,\rho})\leq (E+\epsilon)T$. Notice that here we use that $\lambda$ is at least $C^2$ (so that $X_\lambda$ is at least $C^1$) in order to be able to use Arzela-Ascoli's theorem.

Because $\epsilon>0$ above can be taken arbitrarily close to $0$ we can actually obtain a sequence $\gamma_{j,\rho}$ of Reeb orbits of $X_\lambda$ whose homotopy class is $\rho$ such that the actions $A(\gamma_{j,\rho})$ converges to $ET$. Again applying Arzela-Ascoli's theorem, we obtain that the sequence $\gamma_{j,\rho}$ has a convergent subsequent, which converges to an orbit $\gamma_{\rho}$ satisfying $A(\gamma_{j,\rho}) \leq ET$.

Reasoning as in step 3 above, we obtain that $N_{X_{\lambda}}(T) \geq e^{\frac{aT}{E} + b}$. Applying Theorem \ref{theorem1'} we obtain the desired estimate for the topological entropy. This finishes the proof of the theorem.
\qed

\section{Contact 3-manifolds with a hyperbolic component} \label{section5}

In this section we will prove the following theorem:

\begin{theorem} \label{theorem3'}
Let $M$ be a closed connected oriented  3-manifold which can be cut along a nonempty family of incompressible tori into a family $\{M_i, 0 \leq i \leq k\}$ of irreducible manifolds with boundary, such that the component $M_0$ satisfies:
\begin{itemize}
\item{$M_0$ is the mapping torus of a diffeomorphism $h: S \to S$ with pseudo-Anosov monodromy on a surface $S$ with non-empty boundary.}
\end{itemize}
Then $M$ can be given infinitely many non-diffeomorphic contact structures $\xi_k$, such that for each $\xi_k$ there exists a hypertight contact form $\lambda_k$ on $(M,\xi_k)$ which has exponential homotopical growth of cylindrical contact homology.
\end{theorem}

We denote by $S$ a surface with boundary and $\omega$ a symplectic form on $S$. Let $h$ be a  symplectomorphism of $(S,\omega)$ to itself, with pseudo-Anosov monodromy and which is the identity on a neighbourhood of $\partial S$. We follow a well known recipe to construct a suitable contact form in the mapping torus $\Sigma(S,h)$.

We choose a primitive $\beta$ for $\omega$ such that for coordinates $(r,\theta) \in [-\epsilon,0] \times S^1$ in a neighbourhood $V$ of $\partial S$ we have $ \beta = f(r)d\theta$, where $f>0$ and $f'>0$. We pick a smooth non-decreasing function $ F_0:\mathbb{R} \to [0,1]$   which satisfies $F_0(t)=0$ for $t\in (-\infty,\frac{1}{100})$ and $F_0(t)=1$ for $t\in (\frac{1}{100},+\infty)$. For $i \in \mathbb{Z}$ we define $F_i(t) = F_0(t-i)$.  Fixing $\epsilon >0$, we define a 1-form $\widetilde{\alpha}$ on $\mathbb{R} \times S$ by letting
\begin{equation}
\widetilde{\alpha} = dt + \epsilon(1-F_i(t))(h^{i})^*\beta + \epsilon F_i(t)(h^{i+1})^*\beta \ \mbox{ for } \ t \in [i,i+1)
\end{equation}
It is immediate to see that this defines a smooth 1-form on $\mathbb{R} \times S$, and a simple computation shows that if $\epsilon$ is small enough the 1-form $\widetilde{\alpha}$ is a contact form. For $t\in[0,1]$, the Reeb vector field $X_{\widetilde{\alpha}}$ is equal to $\partial_t + v(p,t)$, where $v(p,t)$ is the unique vector tangent to $S$ that satisfies $\omega(v(p,t), \cdot )= F'_0(t)\beta - F'_0 h^* \beta$.

Consider on the diffeomorphism $H:\mathbb{R} \times S \to \mathbb{R} \times S$ defined by $H(t,p)=(t-1,h(p))$. The mapping torus $\Sigma(S,h)$ is defined by:

\begin{equation}
\Sigma(S,h):=(\mathbb{R} \times S )/_{(t,p) \sim H(t,p)}
\end{equation}
and we denote by $\pi:\mathbb{R} \times S \to \Sigma(S,h)$ the associated covering map.

Because $H^*\widetilde{\alpha}=\widetilde{\alpha}$, there exists a unique contact form $\alpha$ on $\Sigma(S,h)$ such that $\pi^*\alpha =\widetilde{\alpha}$. Notice that in the neighbourhood $S^1 \times V$ of $\partial \Sigma(S,h)$, $\alpha= d t + \epsilon f(r)d\theta$, which implies that $X_{\alpha}$ is tangent to $\partial \Sigma(S,h)$.

The Reeb vector field $X_{\alpha}$ on $\Sigma(S,h)$ is transverse to the surfaces $\{t\} \times S$ for $t\in \mathbb{R}/_{\mathbb{Z}}$. This implies that $\{0\} \times S$ is a global surface of section for the Reeb flow of $\alpha$, and by our expression of $X_{\widetilde{\alpha}}$ the first return map of the Reeb flow of $\alpha $ is isotopic to $h$.

By doing a sufficiently small perturbation of $\alpha$ supported in the interior of $\Sigma(S,h)$ we can obtain a contact form $\widehat{\alpha}$ satisfying that all Reeb orbits of  $\widehat{\alpha}$ which are not freely homotopic to curves in $\partial \Sigma(S,h)$ are non-degenerate, and such that $\{0\} \times S$ is a global surface of section for the flow of $X_{\widehat{\alpha}}$. Notice that as the perturbation is supported in the interior of $\Sigma(S,h)$, the Reeb flow of $\widehat{\alpha}$ is also tangent to the boundary of  $\Sigma(S,h)$.

\subsection{Contact 3-manifolds containing $(\Sigma(S,h),\widehat{\alpha})$ as a component}

Let $M$ be a closed connected oriented  3-manifold which can be cut along a non-empty family of incompressible tori into a family $\{M_i, 0 \leq i \leq k\}$ of irreducible manifolds with boundary, such that the component $M_0$ is diffeomorphic to $(\Sigma(S,h),\alpha)$.
Then it is possible to construct hypertight contact forms on $M$ which match with $\widehat{\alpha}$ in the component $M_0$. More precisely, we have the following result due to Colin and Honda, and Vaugon:

\begin{proposition*} (\cite{CH,Vaugon})
Let $M$ be a closed connected oriented  3-manifold which can be cut along a non-empty family of incompressible tori into a family $\{M_i, 0 \leq i \leq k\}$ of irreducible manifolds with boundary, such that the component $M_0$ is diffeomorphic to $(\Sigma(S,h),\alpha)$. Then, there exist an infinite family $\{\xi_k, k \in \mathbb{Z}\}$ of non-diffeomorphic contact structures on $M$  such that
\begin{itemize}
\item for each $k \in \mathbb{Z}$ there exists a hypertight contact form $\lambda_k$ on $(M,\xi_k)$ which coincides with $\widehat{\alpha}$ on the component $M_0$.
\end{itemize}
\end{proposition*}

We briefly recall the construction of the contact forms $\lambda_k$, and refer the reader to \cite{CH,Vaugon} for the details. When $i \geq 1$, we apply \cite[Theorem 1.3]{CH} to obtain a hypertight contact form $\alpha_i$ on $M_i$ which is compatible with the orientation of $M_i$, and whose Reeb vector field $X_{\alpha_i}$ is tangent to the boundary of $M_i$. On the special piece $M_0$ we consider the the contact form $\alpha_0$ equal to $\widehat{\alpha}$ constructed in the above.

Let $\{\mathfrak{T}_j | 1 \leq \ j \leq m \}$ be the family of incompressible tori along which we cut $M$ to obtain the pieces $M_i$. Then the contact forms $\alpha_i$ give a hypertight contact form on each component of $M \setminus \bigcup_{j \geq1}^m \mathbb{V}(\mathfrak{T}_j) $, where $\mathbb{V}(\mathfrak{T}_j)$ is a small open neighborhood of $T_j$. This gives a contact form $\widehat{\lambda}$ on $M \setminus \bigcup_{j \geq1}^m \mathbb{V}(\mathfrak{T}_j) $.
 Using an interpolation process (see \cite[section 7]{Vaugon}), one can construct contact forms on the neighborhoods $\overline{\mathbb{V}(\mathfrak{T}_j)}$ which coincide with $\widehat{\lambda}$ on $\partial \overline{\mathbb{V}(\mathfrak{T}_j)}$. The interpolation process is not unique and can be done in such ways as to produce  an infinite family of distinct contact forms $\{ \lambda_k \ | \ k\in \mathbb{Z}\}$ on $M$ that extend $\widehat{\lambda}$,  and which are associated to contact structures $\xi_k:= \ker \lambda_k$ that are all non-diffeomorphic. The contact topological invariant used to show that the family $\{\xi_k \ | \ k\in \mathbb{Z}\}$ is composed by non-diffeomorphic contact structures is the \textit{Giroux torsion} (see \cite[section 7]{Vaugon}).

\subsection{Proof of Theorem 3}

It is clear that Theorem \ref{theorem3'} will follow from Theorem \ref{theorem2'}, if we establish that  the  cylindrical contact homology of $\lambda_k$ has exponential homotopical growth. This is the content of

\begin{proposition} \label{proposition2}
$\lambda_k$ has exponential homotopical growth of cylindrical contact homology.
\end{proposition}

Before proving the proposition we introduce some necessary ideas and notation.
The first return map of $X_{\widehat{\alpha}}$ is a diffeomorphism $\widehat{h}:S \to S$ which is homotopic to $h$ and therefore to a pseudo-Anosov map. The Reeb orbits of $X_{\widehat{\alpha}}$ are in one-to-one correspondence with periodic orbits of $\widehat{h}$.  Moreover we have that two Reeb orbits $\gamma_1$ and $\gamma_2$ of $X_{\widehat{\alpha}}$ are freely homotopic if and only if their associated periodic orbits are in the same Nielsen class. Thus there is an injective map $\Xi$ from the set $\mathcal{N}$ of Nielsen classes to the set $\bigwedge$ of free homotopy classes of Reeb orbits in $\Sigma(S,h)$.

Denote by $\mathcal{N}_k$ the set of distinct Nielsen classes which contain only periodic orbits of $\widehat{h}$ of period smaller or equal to $k$. Because of the pseudo-Anosov monodromy of $\widehat{h}$ we know that there are constants $a >0$ and $b\in \mathbb{R}$, such that $\sharp(\mathcal{N}_k) > e^{ak + b}$ for all $k \geq 1$.
Analogously define the subset $\bigwedge_T(\Sigma(S,h))$ of free homotopy classes of $\Sigma(S,h)$ which contains at least one Reeb orbit of $X_{\widehat{\alpha}}$ and contains only Reeb orbits with action smaller than $T$.

Because $\widehat{h}$ is the first return map for a global surface of section of the flow $X_{\widehat{\alpha}}$, there exists a constant $\eta\geq 1$ such that  $\varrho \in \mathcal{N}_k \Rightarrow \Xi(\varrho) \in \bigwedge_{\eta k}(\Sigma(S,h))$. This implies that $\sharp(\bigwedge_T(\Sigma(S,h))) > e^{ \frac{a}{\eta}T+ b} $ for $T\geq \eta$.
Let $\bigwedge^{0}_{T}(\Sigma(S,h))$ be the subset of $\bigwedge_T(\Sigma(S,h))$ which contains free homotopy classes in $\Sigma(S,h)$ which are primitive and different from the ones generated by curves in $\partial\Sigma(S,h)$ (we denote by $\bigwedge^0(\Sigma(S,h))$ the set $\bigwedge^{0}_{+\infty}(\Sigma(S,h))$). Because the fundamental group of $\partial\Sigma(S,h)$ grows quadratically we know that there is $T_0\geq 0$ such that $\sharp{\bigwedge^{0}_{T}}(\Sigma(S,h)) \geq e^{ \frac{a}{\eta}T+ b} $ for $T\geq T_0$..

We are now ready for the proof of Proposition \ref{proposition2}. The main ideas of the argument are due to Vaugon, which estimated in \cite{Vaugon} a different growth rate of the cylindrical contact homology $\lambda_k$.

\

\textit{Proof of Proposition \ref{proposition2}:}

\textbf{Step 1:}

Let $i:\Sigma(S,h) \to M$ be the injection we obtain from looking at $\Sigma(S,h)$ as a component of $M$. Because of the incompressibility of $\partial \Sigma(S,h)$ in $M$, the associated map $i_* : \bigwedge^0_T(\Sigma(S,h)) \to \bigwedge(M)$ is injective for any $T>0$ (here $\Lambda(M)$ denotes the free loop space of $M$). It is clear that all curves belonging to a free homotopy class $\rho \in i_*(\bigwedge^0(\Sigma(S,h)) $ must intersect the component $M_0$.

Using then that the Reeb flow of $\lambda_k$ is tangent to $\partial \Sigma(S,h)$, we conclude that for every $\rho \in i_*(\bigwedge^0(\Sigma(S,h))$ all the Reeb orbits of $X_{\lambda_k}$ that belong to $\rho$ are contained in the interior of $\Sigma(S,h)$. Therefore, the image $i_*(\bigwedge^0(\Sigma(S,h))$ is contained in the et $\bigwedge_{T}(M)$ of free homotopy classes of $M$ which only contain Reeb orbits with action smaller than $T$.

 This mean that the map $i_*:\bigwedge^{0}_{T}(\Sigma(S,h)) \to \bigwedge(M)$, restricts to a map $i_*:\bigwedge^{0}_{T}(\Sigma(S,h)) \to \bigwedge_T(M)$, where by $\bigwedge_{T}(M)$ we denote the set of free homotopy classes of $M$ which only contain Reeb orbits with action smaller than $T$.

\

\textbf{Step 2:} For every $\rho \in i_*(\bigwedge^0(\Sigma(S,h))$ we have $C\mathbb{H}^{\varrho}_{cyl}(\lambda_k) \neq 0$.

Vaugon showed (see the proofs of Lemma 7.11 and Theorems 1.3 and 1.2 in \cite{Vaugon}) that the numbers of even and odd Reeb orbits in $\rho$ differ. For Euler characteristic reasons this implies that $C\mathbb{H}^{\varrho}_{cyl}(\lambda_k) \neq 0$.

\

\textbf{Step 3:}

Recall that in section \ref{section4} we defined $N^{cyl}_T(\lambda_k)$ to be the number of different free homotopy classes $\varrho$ in $\bigwedge_T(M)$ which contained only simple Reeb orbits with action smaller then $T$ and such that $C\mathbb{H}^{\varrho}_{cyl}(\lambda_k) \neq 0$.

Combining the first two steps we obtain
\begin{equation}
N^{cyl}_T(\lambda_k) \geq  \sharp (i_*(\bigwedge^0_T(\Sigma(S,h)))) = \sharp (\bigwedge^0_T(\Sigma(S,h))) \geq e^{\frac{aT}{\eta} + b},
\end{equation}
which establishes the proposition.
\qed

\

\textit{Proof of Theorem \ref{theorem3'}:}
As mentioned previously, Theorem \ref{theorem3'} follows directly from combining Proposition \ref{proposition2} and Theorem \ref{theorem2'}.
\qed

\

It would be interesting to obtain an upper bound on the constant $\eta$ above. This could provide a more precise estimate for the homotopical growth rate of $C\mathbb{H}_{cyl}(\lambda_k)$.

\section{Graph manifolds and Handel-Thurston surgery}\label{section6}

In \cite{HT} Handel and Thurston used Dehn surgery to obtain non-algebraic Anosov flows in 3-manifolds. Their surgery was adapted to the contact setting by Foulon and Hasselblatt in \cite{FH}, who interpreted it as a Legendrian surgery and used it to produce non-algebraic Anosov Reeb flows on 3-manifolds. In this section we apply the Foulon-Hasselblatt Legendrian surgery to obtain more examples of contact 3-manifolds which are distinct from unit tangent bundles, and on which every Reeb flow has positive topological entropy.

Some clarifications regarding the surgeries we consider are in order.  On one hand, we restrict our attention to the Foulon-Hasselblatt surgery on Legendrian lifts of embedded separating geodesics on hyperbolic surfaces. This is an important restriction, since Foulon and Hasselblatt perform their surgery on the Legendrian lift of any immersed closed geodesic on a hyperbolic surface.
On the other hand, for this restricted class of Legendrian knots, the surgery we consider is a bit more general than the one in \cite{FH}. They restrict their attention to Dehn surgeries with positive integer coefficients while we consider the case of any integer coefficient, as is explained in subsection \ref{sec7.1}.

\subsection{The surgery}\label{sec7.1}

We start by fixing some notation. Let $(S,g)$ be an oriented hyperbolic surface and $\mathfrak{r}:S^1 \to S$ an embedded oriented separating geodesic of $g$. We denote by $\pi:(\mathbb{D},g) \to (S,g)$ a locally isometric covering of $(S,g)$ by the the hyperbolic disc $(\mathbb{D},g)$ with the property that $(-1,1) \times \{0\} \subset \pi^{-1}(\mathfrak{r}(S^1))$. Such a covering always exists since the segment $(-1,1) \times \{0\}$ of the real axis is a geodesic in $(\mathbb{D},g)$.
We denote by $v(\theta)$ the unique unitary vector field over $\mathfrak{r}(\theta)$ satisfying $\angle(\mathfrak{r}'(\theta),v(\theta))= -\frac{\pi}{2} $. Our orientation convention is chosen so that for coordinates $z=x+iy\in \mathbb{D}$, the lift of $v(\theta)$ to $(-1,1) \times \{0\} $ is a positive multiple of the vector field $-\partial_y$ over $(-1,1) \times \{0\}$. Also, let $\Pi: T_1 S \to S$ denote the base point projection.

Because  $\mathfrak{r}$ is a separating geodesic, we can cut $S$ along  $\mathfrak{r}$ to obtain two oriented hyperbolic surfaces with boundary which we denote by $S_1$ and $S_2$. Our labelling is chosen so that the vector field $v(\theta)$ points inward $S_2$ and outward $S_1$. This decomposition of $S$ induces a decomposition of $T_1 S $ in $T_1 S_1 $ and $T_1 S_2$. Both  $T_1 S_1 $ and $T_1 S_2$ are 3-manifolds whose boundary is the torus formed by the the unit fibers over $\mathfrak{r}$.

Denote by $V_{\mathfrak{r},\delta}$ the closed $\delta-$neighbourhood of the the geodesic $\mathfrak{r}$ for the hyperbolic metric $g$. For $\delta>0$ sufficiently small we have that $V_{\mathfrak{r},\delta}$ is an annulus such that the only closed geodesics contained in $V_{\mathfrak{r},\delta}$ are the covers of $\mathfrak{r}$, and that satisfies the following convexity property: if $\breve{V}$ is the connected component of $\pi^{-1}(V_{\mathfrak{r},\delta})$ containing $(-1,1) \times \{0\}$, then every segment of a hyperbolic geodesic starting and ending in $\breve{V}$ is completely contained in $\breve{V}$. It also follows from the conventions adopted above, that if we denote by $U^+$ the upper hemisphere of the $\mathbb{D}$ composed of points with positive imaginary component and by $U^+$ the lower hemisphere of the $\mathbb{D}$ composed of points with negative imaginary component, we have:
\begin{equation}
\breve{V}\cap U^+ \subset \pi^{-1}(S_1) \ \mbox{ and } \ \breve{V}\cap U^- \subset \pi^{-1}(S_2).
\end{equation}
This fact has the following important consequence: if $\nu([0,K])$ is a hyperbolic geodesic segment starting and ending at $V_{\mathfrak{r},\delta}$ and contained in one of the $ S_i$, then $[\nu]$ is a non-trivial homotopy class in the relative fundamental group $\pi_1(S_i,V_{\mathfrak{r},\delta})$.

On the unit tangent bundle $T_1 S$ we consider consider the contact form $\lambda_g$ whose Reeb vector field is the geodesic vector field for the hyperbolic metric $g$. It is well known that the lifted curve $(\mathfrak{r}(\theta),v(\theta))$ in $T_1 S$  is Legendrian on the contact manifold $(T_1 S,\ker \lambda_g)$. The  geodesic vector field $X_{\lambda_g}$ over the Legendrian curve coincides with the horizontal lift of $v$ (see \cite[section 1.3]{Pat}), points inward $T_1 S_2$ and outward  $T_1 S_1$, and is normal to $\partial T_1 S_2 = \partial T_1 S_1$ for the Sasaki metric on $ T_1 S$.

Moreover if $\delta >0$ is small enough we know that for every $\vartheta \in L_{\mathfrak{r}}$ there exists numbers $t_1<0$ and $t_2>0$ such that:
\begin{eqnarray}
\phi^{t_1}_{\lambda_g}(\vartheta) \in T_1 S_1 \setminus \Pi^{-1}(V_{\mathfrak{r},\delta}), \\
\phi^{t_2}_{\lambda_g}(\vartheta) \in T_1 S_2 \setminus \Pi^{-1}(V_{\mathfrak{r},\delta}).
\end{eqnarray}

Following \cite{FH}, we know that there exists a neighbourhood $B^{3\eta}_{2\epsilon} $ of $L_{\mathfrak{r}}$ on which we can find coordinates $(t,s,w) \in (-3\eta,3\eta) \times S^1 \times (-2\epsilon,2\epsilon)$ such that:
\begin{eqnarray}
\lambda_g =dt +wds, \\
L_\mathfrak{r} = \{0\} \times S^1 \times \{0\},
\end{eqnarray}
where $\{0\} \times \{\theta\} \times (-2\epsilon,2\epsilon)$ is a local parametrization of the unitary fiber over $\theta \in L_{\mathfrak{r}}$, and $\epsilon < \frac{\eta}{4|q|\pi}$, with $q$ being a fixed integer. Let $\mathcal{W}^-= \{-3\eta\} \times S^1 \times (-2\epsilon,2\epsilon)$ and $\mathcal{W}^+= \{+3\eta\} \times S^1 \times (-2\epsilon,2\epsilon)$. It is clear that $\Pi(\mathcal{W}^-) \subset S_1$ and $\Pi(\mathcal{W}^+) \subset S_2$. Because on $\overline{B}^{3\eta}_{2\epsilon} $ the Reeb vector field $X_{\lambda_g}$ is given by $\partial_t$, it is clear that for every point $p\in B^{3\eta}_{2\epsilon} $ there are $p^- \in \mathcal{W}^-$, $p^+ \in \mathcal{W}^+$, $t^- \in (-6\eta,0)$ and $t^+ \in (0,6\eta)$ for which:
\begin{equation}
\phi^{t^-}_{X_{\lambda_g}}(p) = p^- \ \mbox{and} \ \phi^{t^+}_{X_{\lambda_g}}(p) = p^+
\end{equation}
This means that trajectories of the flow of $X_{\lambda_g}$ that enter the box $B^{3\eta}_{2\epsilon} $ enter through $\mathcal{W}^-$ and exit through $\mathcal{W}^+$. They cannot stay inside $B^{3\eta}_{2\epsilon} $ for a very long positive or negative interval of time. We can say even more about these trajectories.

For $\sigma=(\mathfrak{p},\dot{\mathfrak{p}}) \in (S\times T_p S)$ in $\mathcal{W}^+\cup \mathcal{W}^-$ let $\widetilde{\sigma} = (\widetilde{\mathfrak{p}},\dot{\widetilde{\mathfrak{p}}})$ be a lift of $\sigma$ to the unit tangent bundle $T_1 \mathbb{D}$ such that $\widetilde{\mathfrak{p}} \in \breve{V}$. The geodesic vector field $X_{\lambda_g}$ in $\widetilde{\sigma}$ coincides with the horizontal lift of $\dot{\mathfrak{p}}$ (\cite[section 1.3]{Pat}). For $\delta$, $\eta>0$ and $\epsilon < \frac{\eta}{4|q|\pi}$ sufficiently small we can guarantee that:
\begin{itemize}
\item{$\Pi(B^{3\eta}_{2\epsilon}) $ is contained in $V_{\mathfrak{r},\delta}$,}
\item{for the lifts $\widetilde{\sigma} =(\widetilde{\mathfrak{p}},\dot{\widetilde{\mathfrak{p}}})$ of points in $\mathcal{W}^+\cup \mathcal{W}^-$ as above, the vector $\dot{\widetilde{\mathfrak{p}}}$ (which is the projection of the geodesic vector field $X_{\lambda_g}(\widetilde{\sigma})$) satisfies $\angle(\dot{\widetilde{\mathfrak{p}}}, -\partial_{y}) < \delta$.}
\end{itemize}
With a such a choice of $\delta>0$, $\eta>0$ and $0<\epsilon < \frac{\eta}{4|q|\pi}$, we obtain  that for every $\sigma^+ \in \mathcal{W}^+$ there exists $t_{\sigma^+}>0$ and for every $\sigma^- \in \mathcal{W}^-$ there exists $t_{\sigma^-}<0$  such that:
\begin{eqnarray}
\phi^{t_{\sigma^+}}_{X_{\lambda_g}}(\sigma^+) \in (T_1 S_2) \setminus V_{\mathfrak{r},\delta} \ \mbox{and} \ \forall t\in[0,t_{\sigma^+}] \ \phi^{t}_{X_{\lambda_g}}(\sigma^+) \notin B^{3\eta}_{2\epsilon}, \\ \phi^{t_{\sigma^-}}_{X_{\lambda_g}}(\sigma^-) \in (T_1 S_1) \setminus V_{\mathfrak{r},\delta} \ \mbox{and} \ \forall t\in[t_{\sigma^-},0] \ \phi^{t}_{X_{\lambda_g}}(\sigma^+) \notin B^{3\eta}_{2\epsilon}.
\end{eqnarray}
To prove this last condition above one uses the fact that $\angle(\dot{\widetilde{p}}, -\partial_{y}) < \delta$ is small and studies the behavior of geodesics in $(\mathbb{D},g)$ starting at points close to the real axis and with initial velocity close to $-\partial_y$. It is easy to see that such geodesics have to cut through the region $V_{\mathfrak{r},\delta}$ and visit the interior of both $S_{1} \setminus V_{\mathfrak{r},\delta}$ and $S_{2} \setminus V_{\mathfrak{r},\delta}$.
 From now on we will assume that $\delta>0$, $\eta>0$ and $0<\epsilon < \frac{\eta}{4|q|\pi}$  are such that the all the above mentioned properties described for them being sufficiently small, hold simultaneously.

Consider the map $F:B^{2\eta}_{2\epsilon} \setminus \overline{B}^{\eta}_{\epsilon} \to  B^{2\eta}_{2\epsilon} \setminus \overline{B}^{\eta}_{\epsilon}$ defined by
\begin{equation}
F(t,s,w) = (t,s+f(w),w) \ \mbox{for} \ (t,s,w) \in (\eta,2\eta) \times S^1 \times (-2\epsilon, 2\epsilon),
\end{equation}
where $f(w) = -q\mathcal{R}(\frac{w}{\epsilon})$ (for our previously chosen integer $q$) and  $\mathcal{R}:[-1,1] \to [0,2\pi]$ satisfies $\mathcal{R}=0$ on a neighbourhood of $-1$, $\mathcal{R}=2\pi$ on a neighbourhood of $1$, $0 \leq \mathcal{R}' \leq 4$ and $\mathcal{R}'$ is an even function.

Our new 3-manifold $M$ is obtained by gluing $T_1 S \setminus  \overline{B}^{\eta}_{\epsilon} $ and $B^{2\eta}_{2\epsilon}$ using the map $F$:
\begin{equation}
M= \bigslant{(T_1 S \setminus \overline{B}^{\eta}_{\epsilon}) \cup B^{2\eta}_{2\epsilon} }{(x \in B^{2\eta}_{2\epsilon} \setminus \overline{B}^{\eta}_{\epsilon} ) \thicksim (F(x) \in  T_1 S \setminus  \overline{B}^{\eta}_{\epsilon} )}
\end{equation}

Notice that $T_1 S= \bigslant{(T_1 S \setminus \overline{B}^{\eta}_{\epsilon}) \cup B^{2\eta}_{2\epsilon} }{(x \in B^{2\eta}_{2\epsilon} \setminus \overline{B}^{\eta}_{\epsilon} ) \thicksim (x \in  T_1 S \setminus  \overline{B}^{\eta}_{\epsilon} )}$. This clarifies our construction of $M$ and shows that $M$ is obtained from $T_1 S$ via a Dehn surgery on $L_{\mathfrak{r}}$.
We follow \cite{FH} to endow $M$ with a contact form which coincides $\lambda_g$ outside of $B^{2\eta}_{2\epsilon}$. As a preparation we define the function $\beta: (-3\eta,3\eta) \to \mathbb{R}$:
\begin{itemize}
\item{$\beta$ is equal to $1$ in an open neighbourhood of $[-2\eta,2\eta]$,}
\item{$|\beta'| \leq \frac{\pi}{\eta}$ and $supp \beta $ is contained in $[-3\eta,3\eta]$.}
\end{itemize}

Using $\beta$ we define
\begin{equation}
r(t,w)= \beta(t)\int_{-2\epsilon}^w xf'(x)dx.
\end{equation}
We point out to the reader that $supp(r)$ is contained in $B^{3\eta}_{\epsilon} $ and therefore so is $supp(dr)$. Notice also, that in $B^{2\eta}_{2\epsilon} \setminus \overline{B}^{\eta}_{\epsilon}$ one has $dr = \frac{w}{2}f'(w) dw$.

Again following \cite{FH} we define in $T_1 S \setminus \overline{B}^{\eta}_{\epsilon}$ the 1-form
\begin{eqnarray}
A_r = dt +wds +dr \ \mbox{for} \ (-3\eta,-\eta), \\
A_r = dt + wds -dr \ \mbox{for} \ (\eta,3\eta), \\
A_r = \lambda_g \ \mbox{otherwise}.
\end{eqnarray}
Notice that because $supp(dr)$ is contained in $B^{3\eta}_{\epsilon} $ the 1-form $A_r$ is well-defined.

On the box $B^{2\eta}_{2\epsilon}$ we define
\begin{equation}
\widetilde{A}= dt + wds +dr.
\end{equation}

A direct computation shows that $F^*(A_r) = \widetilde{A}$, which means that the gluing map $F$ allows us to glue the 1-forms $A_r$ and $\widetilde{A}$. We denote by $\lambda_{FH}$ the 1-form in $M$ obtained by gluing $\widetilde{A}$ and $A_r$. We will denote by $\widetilde{B}$ the following region:
\begin{equation}
\widetilde{B}= \bigslant{((B^{3\eta}_{2\epsilon}  \setminus \overline{B}^{\eta}_{\epsilon}) \ \subset M) \cup B^{2\eta}_{2\epsilon} }{(x \in B^{2\eta}_{2\epsilon} \setminus \overline{B}^{\eta}_{\epsilon} ) \thicksim (F(x) \in (B^{3\eta}_{2\epsilon}  \setminus \overline{B}^{\eta}_{\epsilon} )}
\end{equation}
The importance of this region lies in the fact that in $M \setminus \widetilde{B} = T_1 S \setminus B^{3\eta}_{2\epsilon}$, the contact form $\lambda_{FH}$ coincides with $\lambda_g$.

Following \cite{FH} one shows through a direct computation that $(dt + wds \pm dr)\wedge (dw\wedge ds) = (1 \pm \frac{\partial r}{\partial t})dt\wedge dw\wedge ds$. Using the fact that $\epsilon < \frac{\eta}{8\pi |q|}$ one gets that $|\frac{\partial r}{\partial t}|<1$, thus obtaining that $(dt + wds \pm dr)$ is a contact form. It follows from this that $A_r$ and $\widetilde{A}$ are contact forms in their respective domains and therefore $\lambda_{FH}$ is a contact form in $M$. More strongly, Foulon and Hasselblatt proceed to show that if $q$ is non-negative the Reeb flow of $\lambda_{FH}$ is an Anosov Reeb flow.

\subsection{Hypertightness and exponential homotopical growth of contact homology of $\lambda_{FH}$}

For $q \in \mathbb{N}$ the hypertightness of $\lambda_{FH}$ follows from the fact that its Reeb flow is Anosov \cite{Fe}. In this subsection we give an independent and completely geometrical proof of hypertightness of $\lambda_{FH}$, which is valid for every $q \in \mathbb{Z}$.

To understand the topology of Reeb orbits of $\lambda_{FH}$ we will study trajectories that enter the surgery region $\widetilde{B}$.  We start by studying trajectories in $B^{2\eta}_{2\epsilon}$. In this region we have
\begin{equation}
X_{\lambda_{FH}} = \frac{\partial_t}{1 + \partial_t r}.
\end{equation}
This implies, similarly to what happens of $\lambda_g$, that for points $p \in B^{2\eta}_{2\epsilon}$ the trajectory $\phi^{t}_{X_{\lambda_{FH}}}(p)$ leaves the box $B^{2\eta}_{2\epsilon}$ in forward and backward times. More precisely, there exists a constant $\widetilde{a}>0$ depending only on $\lambda_{FH}$, such that for $p \in B^{2\eta}_{2\epsilon}$ there are $\breve{p}^- \in \breve{\mathcal{W}}^-= \{-2\eta\} \times S^1 \times [-2\epsilon,2\epsilon]$, $\breve{p}^+ \in \breve{\mathcal{W}}^+= \{+2\eta\} \times S^1 \times [-2\epsilon,2\epsilon]$, $\breve{t}^- \in (-\widetilde{a},0]$ and $\breve{t}^+ \in [0,\widetilde{a})$ such that
\begin{equation*}
\phi^{t}_{X_{\lambda_{FH}}}(\breve{p}) \mbox{ is in the interior of } B^{2\eta}_{2\epsilon} \mbox{ for every } t \in (t^-,t^+),
\end{equation*}
\begin{equation}
\phi^{t^-}_{X_{\lambda_{FH}}}(\breve{p}) = \breve{p}^- \ \mbox{and} \ \phi^{t^+}_{X_{\lambda_g}}(\breve{p}) = \breve{p}^+.
\end{equation}

We now analyse the trajectories of points $\breve{p}^- \in \breve{\mathcal{W}}^-$ and $\breve{p}^+ \in \breve{\mathcal{W}}^+$. For this, we first notice that on $\widetilde{B} \setminus B^{\eta}_{\epsilon}$ the contact form $\lambda_{FH}$ is given by
$dt+wds \pm dr$, and therefore we have in this region
\begin{equation}
X_{\lambda_{FH}}=\frac{\partial_t}{1 \pm \partial_t r},
\end{equation}
which is still a positive multiple of $\partial_t$.

This implies that for every $\breve{p}^- \in \breve{\mathcal{W}}^-$ and $\breve{p}^+ \in \breve{\mathcal{W}}^+$ there exist $t^{\breve{p}^- }<0$ and $t^{\breve{p}^+}<0$ such that
\begin{equation}
\phi^{t^{\breve{p}^- }}_{X_{\lambda_{FH}}}(\breve{p}^-) \in \mathcal{W}^- \ \mbox{and} \ \phi^{t^{\breve{p}^- }}_{X_{\lambda_g}}(\breve{p}^+) \in \mathcal{W}^+
\end{equation}
Again using that $X_{\lambda_{FH}}$ is a positive multiple of $\partial_t$ on $\widetilde{B} \setminus B^{2\eta}_{2\epsilon}$
we have that for every point $p$ in $\widetilde{B} \setminus B^{2\eta}_{2\epsilon}$ whose $t$ coordinate is in $[2\eta,3\eta]$ the trajectory of the flow $\phi^t_{X_{\lambda_{FH}}}$ going through $p$ is a straight line with fixed coordinates $s$ and $w$, that goes from $\breve{\mathcal{W}}^+$ to $\mathcal{W}^+$. Analogously, for every point $p$ in $\widetilde{B} \setminus B^{2\eta}_{2\epsilon}$ whose $t$ coordinate is in $[-3\eta,-2\eta]$ the trajectory of the backward flow of $\phi^t_{X_{\lambda_{FH}}}$ going through $p$ is a straight line $\breve{\mathcal{W}}^-$ to $\mathcal{W}^-$.

Summing up, with all the cases considered above we have showed that for every point $p \in \widetilde{B}$ the trajectory of the flow $\phi^t_{X_{\lambda_{FH}}}$ going through $p$ for $t=0$ intersects $\mathcal{W}^-$ for non-positive time and $\mathcal{W}^+$ for for non-negative time. In other words, all trajectories that intersect $\widetilde{B}$ enter through $\mathcal{W}^-$ and leave through $\mathcal{W}^+$, which means that for all $\check{p} \in \widetilde{B}$ there exist times $t^-_{\check{p}}\leq 0$ and $t^+_{\check{p}}\geq 0$ such that
\begin{eqnarray}
\phi^{t^+_{\check{p}}}_{X_{\lambda_{FH}}}(\check{p}) \in  \mathcal{W}^+, \  \  \  \  \  \  \   \  \  \  \  \  \  \  \  \  \ \\
\phi^{t^-_{\check{p}}}_{X_{\lambda_{FH}}}(\check{p}) \in  \mathcal{W}^-,   \  \  \  \  \  \  \   \  \  \  \  \  \  \  \  \  \   \\
\phi^{t}_{X_{\lambda_{FH}}}(\check{p}) \in \widetilde{B} \mbox{ for all } t \in [t^-_{\check{p}},t^+_{\check{p}}].
\end{eqnarray}

Now, because on $M\setminus \widetilde{B} = T_1 S \setminus  B^{3\eta}_{2\epsilon}$ the contact form $\lambda_{FH}$ coincides with $\lambda_g$ we have that trajectories of $X_{\lambda_{FH}}$ starting at $\mathcal{W}^-$ at the time $t=0$ have to leave $M \setminus N$ as time diminishes before reentering on $\widetilde{B}$. Similarly the trajectories starting at $\mathcal{W}^+$ have to leave $M \setminus N$ for positive time before reentering on $\widetilde{B}$. More precisely, one can use equations (51) and (52) to show that for $p^- \in \mathcal{W}^-$ and $p^+ \in \mathcal{W}^+$ there exist $t_{p^-}<0$ and $t_{p^+}>0$ such that
\begin{eqnarray}
\phi^{t_{p^+}}_{X_{\lambda_{FH}}}(p^+) \in M_2 \setminus N \ \mbox{and} \ \forall t\in[0,t_{p^+}] \ \phi^{t}_{X_{\lambda_{FH}}}(p^+) \notin \widetilde{B}, \\ \phi^{t_{p^-}}_{X_{\lambda_{FH}}}(p^-) \in M_1 \setminus N \ \mbox{and} \ \forall t\in[t_{p^-},0] \ \phi^{t}_{X_{\lambda_{FH}}}(p^-) \notin \widetilde{B},
\end{eqnarray}
where
\begin{eqnarray}
M_1 = \bigslant{(T_1 S_1 \setminus B^{\eta}_{\epsilon}) \cup B^{2\eta}_{2\epsilon}(-)}{(x \in B^{2\eta}_{2\epsilon}(-) \setminus \overline{B}^{\eta}_{\epsilon} ) \thicksim (F(x) \in ((B^{3\eta}_{2\epsilon}\cap T_1 S_1)  \setminus \overline{B}^{\eta}_{\epsilon} )}, \\ M_2 = \bigslant{(T_1 S_2 \setminus B^{\eta}_{\epsilon}) \cup B^{2\eta}_{2\epsilon}(+)}{(x \in B^{2\eta}_{2\epsilon}(+) \setminus \overline{B}^{\eta}_{\epsilon} ) \thicksim (F(x) \in ((B^{2\eta}_{2\epsilon}\cap T_1 S_2)  \setminus \overline{B}^{\eta}_{\epsilon} )}, \\ N= \bigslant{\Pi^{-1}(V_{c,\delta} \setminus) \cup B^{2\eta}_{2\epsilon}(-)}{(x \in B^{2\eta}_{2\epsilon}(-) \setminus \overline{B}^{\eta}_{\epsilon} ) \thicksim (F(x) \in ((B^{3\eta}_{2\epsilon}\cap T_1 S^1)  \setminus \overline{B}^{\eta}_{\epsilon} )},
\end{eqnarray}
for $B^{2\eta}_{2\epsilon}(-) = [-2\eta,0]\times S^1 \times (-2\epsilon,2\epsilon)$ and $B^{2\eta}_{2\epsilon}(+) = [0,2\eta]\times S^1 \times (-2\epsilon,2\epsilon)$.

\textbf{Remark:} It is not hard to see that $M=\bigslant{M_1 \cup M_2}{(x \in \partial M_1) \thicksim (\widetilde{F}(x) \in \partial M_2)}$. Here $\widetilde{F}$ is a Dehn twist which coincides with $(s+f(w),w)$ for $w\in[-2\epsilon,2\epsilon]$ and is the identity elsewhere. This picture of $M$ is closer to the one in the paper \cite{HT} and shows that $M$ is a graph manifold (a graph manifold is one whose JSJ decomposition consists of Seifert $S^1$ bundles). By using this description of $M$ and applying Van-Kampen's to analyse the fundamental group of $M$, Handel and Thurston show that, for $q$ not belonging to a finite subset of $\mathbb{Z}$, no finite cover of $M$ is a Seifert manifold thus obtaining that $M$ is an ``exotic'' graph manifold.

From their definition one sees that as manifolds, $M_1\cong T_1 S_1$ and $M_2 \cong T_1 S_2$. This implies $\partial M_1$ and $\partial M_2$ are incompressible tori in, respectively, $M_1$ and $M_2$. If we look at $M_1$ and $M_2$ as submanifolds of $M$, their boundary $\mathbb{T}$ coincide and is also incompressible in $M$. We remark that $M_i \setminus N$ is diffeomorphic to $T_1 S_i \setminus \Pi^{-1}(V_{c,\delta})$ which is diffeomorphic to $T_1 S_i$ for $i=1,2$.

In a similar way we can describe the topology of $N$. Let $N_i= M_i \cap N$. Reasoning identically as one does to show that $M_i$ is diffeomorphic to $T_1 S_i$ one shows that $N_i$ is diffeomorphic to a thickned two torus $ \mathcal{T}^2 \times [-1,1]$. As $N$ is obtained from $N_1$ and $N_2$ by gluing them along $\mathbb{T}$ (which is a boundary component of both of them) we have that $N$ is also diffeomorphic to the product $\mathcal{T}^2 \times [-1,1]$ .

The discussion above proves the following
\begin{lemma}
For all $\breve{p} \in \widetilde{B}$ the trajectory $\{\phi^{t}_{X_{\lambda_{FH}}}(\breve{p}) \ | \ t \in \mathbb{R}\}$ intersects $M_1 \setminus N$ and $M_2 \setminus N$.
\end{lemma}
\textit{Proof:} We have already established that for $\breve{p} \in \widetilde{B}$ its trajectory intersect $\mathcal{W}^+$ for some non-negative time and $\mathcal{W}^-$ for some non-positive time, as it is shown in equations (65) and (66). One now applies equations (68) and (69) to finish the proof of the lemma.
\qed

\

Notice that trajectories can only enter in $\widetilde{B}$ through the wall $\mathcal{W}^-$ which is contained in $M_1$ and can only exit $\widetilde{B}$ through the wall $\mathcal{W}^+$ which is contained in $M_2$. We also point out that all trajectories of the flow $\phi^{t}_{X_{\lambda_{FH}}}$ are transversal to $\mathbb{T}$, with the exception of the two Reeb orbits which correspond to parametrizations of the hyperbolic geodesic $\mathfrak{r}$ (they continue to exist as periodic orbits after the surgery because they are distant from the surgery region).

We will deduce from the previous discussion the following important lemma.
\begin{lemma} \label{lemma2}
Let $\gamma([0,T'])$ be a trajectory of $X_{\lambda_{FH}}$ such that $\gamma(0) \in \mathbb{T}$,  $\gamma(T') \in \mathbb{T}$ and for all $t \in (0,T')$ we have $\gamma(t) \notin \mathbb{T}$ (notice that in such a situation $\gamma([0,T']) \subset M_i$ for some $i$ equals to $1$ or $2$). Then $\gamma([0,T']) \cap (M_i \setminus N)$ is non-empty.
\end{lemma}

\textit{Proof:} We divide the proof in 3 possible scenarios.

\textbf{First case:} suppose that $\gamma([0,T'])\cap \widetilde{B}$ is empty. In this case $\gamma([0,T'])$ also exists as a hyperbolic geodesic with endpoints in the closed geodesic $\mathfrak{r}$. It follows from the convexity of the hyperbolic metric that $[\gamma([0,T'])] \in \pi_1(T_1 S_i,\mathbb{T})$ is non-trivial. This implies that $[\gamma([0,T'])] \in \pi_1(M_i,\mathbb{T})$ is non-trivial which can be true only if $\gamma([0,T'])\cap (M_i \setminus N)$ is non-empty since $N$ is a tubular neighbourhood of $\mathbb{T}$.

\textbf{Second case:} suppose that $\gamma([0,T'])\cap \widetilde{B}$ is non-empty and $\gamma([0,T']) \subset M_2$. Take $\widehat{t} \in [0,T']$ such that $\gamma(\widehat{t}) \in \widetilde{B}$. We know from our previous discussion that there are $\widehat{t}_1 \leq \widehat{t} \leq \widehat{t}_2$ such that $\gamma([\widehat{t}_1,\widehat{t}_2]) \subset \widetilde{B} $, $\gamma(\widehat{t}_1) \in (\mathbb{T}\cap \widetilde{B}) $ and $\gamma(\widehat{t}_2) \in \mathcal{W}^+$; notice that in the coordinates $(t,s,w)$ for $\widetilde{B}$ considered previously, $\mathbb{T}\cap \widetilde{B}$ is the annulus $\{0\} \times S^1 \times (-2\epsilon,2\epsilon)$. From this picture it is clear that for $t$ smaller that $\widehat{t}_1$ the trajectory enters in $M_1$. Therefore we must have $\widehat{t}_1 =0 $ and $\gamma([0,\widehat{t}_2]) \subset \widetilde{B} $. Notice also that for all $t$ slightly bigger than $\widehat{t}_2$ the trajectory is outside $\widetilde{B}$. Because trajectories of $X_{\lambda_{FH}}$ can only enter $\widetilde{B}$ in $M_1$ we obtain that $\gamma([\widehat{t}_2,T'])$ does not intersect the interior of $\widetilde{B}$ and therefore exists as a hyperbolic geodesic in $T_1 S_2$. Now, using equations (51) and (52) we obtain that, because $\gamma(\widehat{t}_2) \in \mathcal{W}^+$, the trajectory $\gamma: [\widehat{t}_2,T'] \to M_2$ has to intersect $M_2 \setminus N$ before hitting $\mathbb{T}$ at $t=T'$. Thus there is some $t\in (\widehat{t}_2,T') $ for which $\gamma(t) \in M_2 \setminus N$.

\textbf{Third case:} the proof in the case where $\gamma([0,T'])\cap \widetilde{B} $ is non-empty and $\gamma([0,T']) \subset M_1$ is  analogous to the one of the Second case.

\

This three cases exhaust all possibilities and therefore prove the lemma.
\qed

Our reason for introducing the above decomposition of $M$ into $M_1$ and $M_2$ and for proving the lemmas above is to introduce the following representation of Reeb orbits of $\lambda_{FH}$. Let $(\gamma,T)$ be a Reeb orbit of $\lambda_{FH}$ which intersects both $M_1 \setminus N$ and $M_2 \setminus N$. We can assume that the chosen parametrization of the Reeb orbit is such that $\gamma(0) \in \partial N $, and that there are $t_+>0$ and $t_-<0$ such that:
\begin{eqnarray}
\gamma(t_+) \in M_1 \setminus N \ \mbox{and} \ \gamma([0,t_+]) \in M_1 \cup N, \\ \gamma(t_-) \in M_2 \setminus N \ \mbox{and} \ \gamma([t_-,0]) \in M_2 \cup N.
\end{eqnarray}
This means that in an interval of the origin $\gamma$ is coming from $M_2 \setminus N$ and going to $M_1 \setminus N$. It follows from Lemma \ref{lemma2} that there exists a unique sequence $0=t_0<t_{\frac{1}{2}}< t_1<t_{\frac{3}{2}}< ...<t_n=T$ such that $\forall k \in \{0,...,n-1\}$:
\begin{itemize}
\item{$\gamma([t_k,t_{k+\frac{1}{2}}]) \subset M_i$ for $i$ equals to $1$ or $2$,}
\item{$\gamma([t_{k+\frac{1}{2}},t_{k+1}])\in N$ and there is a unique $\widetilde{t_k} \in [t_{k+\frac{1}{2}},t_{k+1}]$ such that $\gamma(\widetilde{t_k}) \in \mathbb{T}$,}
\item{if $\gamma([t_k,t_{k+\frac{1}{2}}]) \subset M_i$ then $\gamma([t_{k+1},t_{k+\frac{3}{2}}]) \subset M_j$ for $j\neq i$.}
\end{itemize}
Notice that $\gamma([t_0,t_{\frac{1}{2}}]) \subset M_1$ and $\gamma([t_{n-1},t_{n-\frac{1}{2}}]) \subset M_2$. This implies that $n$ is even so that we can write $n = 2n'$, and that $\gamma([t_k,t_{k+\frac{1}{2}}]) \subset M_1$ for $k$ even, and $\gamma([t_k,t_{k+\frac{1}{2}}]) \subset M_2$ for $k$ odd. For each $k \in \{0,...,2n'-1\}$ the existence of the unique $\widetilde{t_k}$ in the interval $[t_{k+ \frac{1}{2}},t_{k+1}]$ for which $\gamma(\widetilde{t_k})\in \mathbb{T}$ is guaranteed from Lemma \ref{lemma2} and the fact that $\mathbb{T}$ is the hypersurface that separates $M_1$ and $M_2$.

In order to obtain information on the free homotopy class of $(\gamma,T)$ we observe that for $\gamma([t_k,t_{k+\frac{1}{2}}])$ coincides with a hyperbolic geodesic segment in $T_1 S_i$ starting and ending $V_{\mathfrak{r},\delta}$. Therefore, as we have previously seen the homotopy class $[\gamma([t_k,t_{k+\frac{1}{2}}])]$
in $\pi_{1}(T_1 S_i,V_{\mathfrak{r},\delta})$ is non-trivial which implies that $\gamma([t_k,t_{k+\frac{1}{2}}])$ is a non-trivial relative homotopy class in $\pi_{1}(M_i,N)$. We consider now the curve $\gamma([\widetilde{t}_k,\widetilde{t}_{k+1}])$: it is the concatenation of 3 curves, the first and the third ones being completely contained in $N$ and the middle one being $\gamma([t_k,t_{k+\frac{1}{2}}])$; from this description and the fact that $\gamma([t_k,t_{k+\frac{1}{2}}])$ is a non-trivial relative homotopy class in $\pi_{1}(M_i,N)$ it is clear that $\gamma([\widetilde{t}_k,\widetilde{t}_{k+1}])$ is also non-trivial in $\pi_{1}(M_i,N)$ (and also non-trivial in $\pi_{1}(M_i,\mathbb{T})$).

We now denote by $\widetilde{M}$ the universal cover of $M$ and and $\widehat{\pi} : \widetilde{M} \to M$ a covering map. From the incompressibility of $\mathbb{T}$ it follows that every lift of $\mathbb{T}$ is an embedded plane in $\widetilde{M}$. We denote by $\widetilde{N}^0$ a lift of $N$. Because $N$ is a thickened neighbourhood of an incompressible torus it follows that $\widetilde{N}^0$ is diffeomorphic to $\mathbb{R}^2 \times [-1,1]$, i.e. it is a thickened neighbourhood of an embedded plane in $\widetilde{M}$. Because $N$ separates $M$ in two components, it follows that $\widetilde{N}^0$ separates $\widetilde{M}$ is two connected components. $\partial \widetilde{N}^0$ is the union of two embedded planes $P^0_-$ and $P^0_+$ which are characterized by the fact that there are neighbourhoods $V_-$ and $V_+$ of, respectively, $P^0_-$ and $P^0_+$ such that $\widehat{\pi}(V_-) \subset M_1$ and $\widehat{\pi}(V_+) \subset M_2$. We will denote by $C^0_-$ the connected component of $\widetilde{M} \setminus \widetilde{N}^0$ which intersects $V_-$, and by $C^0_+$ the connected component of $\widetilde{M} \setminus \widetilde{N}^0$ which intersects $V_+$.

As we saw earlier, $[\gamma([t_k,t_{k+ \frac{1}{2}}])]$ is a non-trivial relative homotopy class in $\pi_1(M_i,N)$. We show that this class remains non-trivial when seen in $\pi_1(M,N)$. Let $\mathbb{T}_i=\partial N \cap M_i$.  Because $N$ is obtained by attaching over each point of
$\mathbb{T}_i$ a small compact interval (i.e it is a bundle over $\mathbb{T}_i$ whose fibers are intervals) it follows that $[\gamma([t_k,t_{k+ \frac{1}{2}}])$ would be trivial in $\pi_1(M_i,\mathbb{T}_i)$ if, and only if, it is trival in $\pi_1(M_i,N)$, which is not the case.  As $\mathbb{T}_i$ is isotopic to $\mathbb{T}$, it is also an incompressible torus that divide $M$ in two components. Now, $[\gamma([t_k,t_{k+ \frac{1}{2}}])]$ would be trivial in $\pi_1((M_i \setminus int(N)),\mathbb{T}_i)$ if, and only if, there existed a curve $\mathfrak{c}$ in $\mathbb{T}_i$ with endpoints $\gamma(t_k)$ and $\gamma(t_{k+ \frac{1}{2}})$, such that the concatenation $\gamma \ast \mathfrak{c}$ was contractible in $(M_i \setminus int(N))$. Because of the incompressibility of $\mathbb{T}_i$ such a curve $\gamma \ast \mathfrak{c}$ should be contractible in $(M_i \setminus int(N))$ if, and only if, it was contractible in $M$. This implies that $[\gamma([t_k,t_{k+ \frac{1}{2}}])]$ would be trivial in $\pi_1(M,\mathbb{T}_i)$ if, and only if, it was trivial in $\pi_1((M_i \setminus int(N)),\mathbb{T}_i)$ which we know not to be the case. Lastly, again because $N$ is an interval bundle over $\mathbb{T}_i$, it is clear that as  $[\gamma([t_k,t_{k+ \frac{1}{2}}])]$ is not trivial in $\pi_1(M,\mathbb{T}_i)$ it cannot be trivial in $\pi_1(M,N)$, as we wished to show.

Let now $\widetilde{\gamma}$ be a lift of $\gamma$ such that $\widetilde{\gamma}(0) \in \widetilde{N}^0$. We know that $\widetilde{\gamma}([t_{2n'-\frac{1}{2}} - T , t_{\frac{1}{2}}]) \subset \widetilde{N}^0$. It will be useful to us to define the following sequence:
\begin{equation}
\widetilde{t}_i = q_i T + t_{r_i},
\end{equation}
where $q_i$ and $r_i<2n'$ are the unique integers such that $i = q_i (2n') + r_i$.
Associated to $\widetilde{t}_i$ we associate the lift $\widetilde{N}^i$ of $N$, which is determined by the property that $\widetilde{\gamma}(\widetilde{t}_i) \in \widetilde{N}^i$. It is clear that the sequence $\widetilde{N}^i$ contains all lifts of $N$ which are intersected by the curve $\widetilde{\gamma}(\mathbb{R})$. For the lifts $\widetilde{N}^i$ we define the connected components $C^i_-$ and $C^i_+$ of $\widetilde{M} \setminus \widetilde{N}^i$, and the planes $P^i_-$ and $P^i_+$ analogously as how we defined them for $\widetilde{N}^0$. A priori it could be that for $i\neq j$ we had $\widetilde{N}^i =\widetilde{N}^j$. We will show however, that this cannot happen.

Firstly, $\widetilde{N}^0 \neq \widetilde{N}^1$because $\gamma([\widetilde{t}_0,\widetilde{t}_1])$ is non-trivial in $\pi_1(M,N)$. Also, we have that $\widetilde{N}^1 \subset C^0_-$ because $\gamma([t_0,t_{\frac{1}{2}}]) \subset M_1$. An identical reasoning shows that $\widetilde{N}^2 \neq \widetilde{N}^1$ and
\begin{equation}
\widetilde{N}^2 \subset C^1_{+}.
\end{equation}
On the other hand we have that $\widetilde{N}^0 \subset C^1_-$, because $\widetilde{\gamma}([\widetilde{t}_0 , t_{\frac{1}{2}} ])$ gives a path totally contained in $\widetilde{M} \setminus \widetilde{N}^1$ connecting $\widetilde{N}^0$ and $P^1_-$. As $\widetilde{N}^2 \subset C^1_+$ and $\widetilde{N}^0 \subset C^1_-$, we must have $\widetilde{N}^2 \neq \widetilde{N}^0$. In an identical way, one shows that $\widetilde{N}^3 \neq \widetilde{N}^1$, and more generally that $\widetilde{N}^{i+2} \neq \widetilde{N}^{i}$ and $\widetilde{N}^{i+1} \neq \widetilde{N}^i$.
Now for $\widetilde{N}^3$, we have that $\widetilde{N}^3 \subset C^2_-$. As $\widetilde{\gamma}([\widetilde{t}_0,t_{\frac{3}{2}}])$ is a path completely contained in $\widetilde{M} \setminus \widetilde{N}^2$ connecting $\widetilde{N}^0 $ and $P^2_+$ we obtain that $\widetilde{N}^0 \subset C^2_+$, and therefore $\widetilde{N}^3 \neq \widetilde{N}^0$.

Proceeding inductively along this line one obtains that $\widetilde{N}^i \neq \widetilde{N}^0$ for all $i \neq 0$, and more generally, $\widetilde{N}^i \neq \widetilde{N}^j $ for all $i\neq j$. As a consequence of this, we obtain that the curve $\widetilde{\gamma}(\mathbb{R})$ cannot be homeomorphic to a circle and therefore $\gamma(\mathbb{R})$ cannot be contractible.
We are ready for the main result of this subsection.
\begin{proposition} \label{proposition3}
$\lambda_{FH}$ is hypertight.
\end{proposition}
\textit{Proof:} there are two possibilities for Reeb orbits.

\

\textbf{Possibility 1:} the Reeb orbit $\gamma$ visits both $M_1 \setminus N$ and $M_2 \setminus N$.

In this case, we have just showed above that $\gamma$ is not contractible.

\

\textbf{Possibility 2:} the Reeb orbit $\gamma$ is completely contained in $M_i$ for $i$ equal to $1$ or $2$.

In this case, the Reeb orbit does not visit the surgery region $\widetilde{B}$. Therefore it existed also before the surgery as a closed hyperbolic geodesic in $M_i \setminus \widetilde{B} = T_1 S_i \setminus B^{3\eta}_{2\epsilon}$. Such a closed geodesic is non-contractible in $T_1 S_i$ which is diffeomorphic to $M_i$. We have thus obtained that $\gamma \subset M_i$ is non-contractible in $M_i$.

Looking now at $M_i$ as a submanifold with boundary of $M$, we recall that $\partial M_i$ is an incompressible torus in $M$. This implies that every non-contractible closed curve in $M_i$ remains non-contractible in $M$. Therefore $\gamma$ is also a non-contractible Reeb orbit for this case.
\qed

\

\subsubsection{Exponential homotopical growth of cylindrical contact homology for $\lambda_{FH}$}

We proceed now to obtain more information on the properties of periodic orbits of $X_{\lambda_{FH}}$. We state the following important fact:

\
\begin{lemma} \label{lemma3}
If a Reeb orbit $(\gamma,T)$ of $\lambda_f$ visits both $M_1 \setminus N$ and $M_2 \setminus N$, then any curve freely homotopic to  $(\gamma,T)$ must always intersect $\mathbb{T}$.
\end{lemma}
\textit{Proof of lemma:} As we saw earlier the lift $\widetilde{\gamma}$ intersects all the elements of the sequence $\widetilde{N}_i $ (of lifts of $N$) which satisfy $\widetilde{N}_i \neq \widetilde{N}_j$ for all $ i \neq j$.

Introducing an auxiliary distance $d$ on the compact manifold $M $ (coming from a Riemannian metric) we obtain an auxiliary distance $\widetilde{d}$ on $\widetilde{M}$ by pulling $d$ back by the covering map. It is clear that for $i$ sufficiently big the $\widetilde{d}$ distance between $\widetilde{N}_{\pm i}$ and $\widetilde{N}_0$ becomes arbitrarily large.
As a consequence, one obtains that for each $K> 0$ there exists $t_K>0$ such that $\widetilde{d}(\widetilde{\gamma}(\pm t_K), \widetilde{N}_0) > K$.

Let now $\zeta: [0,T] \to M$ be closed curve freely homotopic to $\gamma([0,T])$. An homotopy $H:[0,T] \times [0,1] \to M$ generates an homotopy $\widetilde{H}: \mathbb{R} \times [0,1] \to \widetilde{M}$ from a lift $\widetilde{\gamma}$ and a lift $\widetilde{\zeta}$. Using the fact that $H$ is uniformly continuous one proves that there exists a constant $\mathfrak{C}>0$ such that $\widetilde{d}(\widetilde{H}(\{t\} \times [0,1]),\widetilde{\gamma}(t))<\mathfrak{C}$ for all $t \in \mathbb{R}$.

Take now $K> 2\mathfrak{C}$. Using the triangle inequality and the facts that $\widetilde{d}((\widetilde{H}(\{t\} \times [0,1])),\widetilde{\gamma}(t))<\mathfrak{C}$ and $\widetilde{d}(\widetilde{\gamma}(\pm t_K), \widetilde{N}_0) > K$ we obtain $H(\{t_K\}\times [0,1])$ is always in the same connected component of $\widetilde{\gamma}(t_K)$. This implies that $\widetilde{\zeta}(\mathbb{R})$ visits both connected components of $\widetilde{M} \setminus \widetilde{N}_0$ and must thus intersect $ \widetilde{N}_0$. Even more, because $\widetilde{\zeta}(\mathbb{R})$ intersects both components of $\partial \widetilde{N}_0 $ we have that $\zeta$ visits both components of $M \setminus N$ and therefore has to intersect $\mathbb{T}$. This completes the proof of the lemma.
\qed

\

We are now ready for the most important result of this section:
\begin{theorem} \label{theorem4'}
Let $(M,\xi_{(q,c)}))$ be the contact manifold endowed obtained by the Foulon-Hasselblatt surgery, and $\lambda_{FS}$ be the contact form obtained via the Foulon-Hasselblat surgery on the Legendrian lift $L_{\mathfrak{r}} \subset T_1 S$. Then $\lambda_{FH}$ is hypertight and its cylindrical contact homology has exponential homotopical growth.
\end{theorem}
We divide the proof that $C\mathbb{H}^{cyl}(M,\lambda_{FH})$ has exponential homotopical growth in steps.

\

\textbf{Step 1:} A special class of Reeb orbits.

We will obtain our estimate by looking at Reeb orbits which are completely contained in the component $M_1$. As we saw previously such orbits never cross the surgery region $\widetilde{B}$. Thus they are in a region where $\lambda_{FH}$ coincides with $\lambda_g$, and such Reeb orbits exist also as closed geodesics in $(S_1,g)$. Conversely, every closed geodesic in $(S_1,g)$ does not cross the region $B_{2\epsilon}^{3\eta}$ and thus also exist as Reeb orbit of $\lambda_{FH}$. This gives a bijective correspondence between closed geodesics of $(S_1,g)$ which are not homotopic to a multiple of $\partial S_1$ and Reeb orbits of $\lambda_{FH}$ which are completely contained in $M_1$.

Let $\bigwedge(S_1)$ denote the set of free homotopy classes in $S_1$ which are not covers of $[\partial S_1]$. We know that each $\rho \in \bigwedge(S_1)$ contains exactly one closed geodesic $c_\rho$. Letting $\gamma_\rho$ be the canonical lift of $c_\rho$ to $T_1 S_1$, we know that $\gamma_\rho$ is a Reeb orbit of $\lambda_g$. As we saw above each $\gamma_\rho$ can also be seen as a Reeb orbit of $\lambda_{FH}$. We will denote by $\bigwedge(S_1)^{\leq T}$ the set primitive of free homotopy classes in $S_1$ whose unique closed geodesic has period smaller or equal to $T$. Because $g$ is hyperbolic it is a well known fact that there exist constants $a>0$, $b$ such that $\sharp(\bigwedge(S_1)^{\leq T})\geq e^{aT+b}$.

Let $\Theta:\bigwedge(S_1) \to \bigwedge(T_1 S_1)$ (where $\bigwedge(T_1 S_1)$ is the free loop space of $T_1 S_1$), be the map which associates $c_\rho$ to $\gamma_\rho$ in $T_1 S_1$. $\Theta:\bigwedge(S_1) \to \bigwedge(T_1 S_1)$  is easily seen to be injective. Because $T_1 S_1$ is diffeomorphic to $M_1$ we can also view $\Theta(\bigwedge(S_1))$ as a subset of the free loop space $\bigwedge(M_1)$ of $M_1$.

\

\textbf{Step 2:}

Let $i:M_1 \to M$ be the injection obtained by looking at $M_1$ as a component of $M$. As seen before the boundary $\partial(i(M_1))=\mathbb{T}$ is an incompressible torus in $M$. We consider the induced map of free loop spaces $i_*:\bigwedge(M_1) \to \bigwedge(M)$. As a consequence of the incompressibility of $\partial(i(M_1)$, the restriction of $i_*$ to $\Theta(\bigwedge(S_1))$ is injective.

To see that, it suffices to show the following claim: if $\zeta$ and $\zeta'$ are curves in $M_1$ which cannot be isotoped to a curve in $\partial M_1$ and which are in the same free homotopy class in $M$, then $\zeta$ and $\zeta'$ are freely homotopic in $M_1$. For $\zeta$ and $\zeta'$ satisfying the hypothesis of our claim there is a cylinder $\mathrm{Cyl}$ in $M$ whose boundary components are $\zeta$ and $\zeta'$ which intersects $\partial M_1$ transversely. In such a case, $\mathrm{Cyl}$ intersects $\partial M_1$ in a finite collection of curves $\{w_n\}$ which are all contractible in $M$; the contractibility of these curves is due to the fact that both $\zeta$ and $\zeta'$ cannot be isotoped to a curve contained in $\partial M_1$. The incompressibility of $\partial M_1$ implies that these $\{w_n\}$ are all contractible already in $\partial M_1$. Now, we cut the discs in $\mathrm{Cyl}$ whose boundary are the curves $c_n$ and substitute them by discs contained in $\partial M_1$. This produces a cylinder $\mathrm{Cyl}'$ completely contained in $M_1$ whose boundaries are $\zeta$ and $\zeta'$. This implies that $\zeta$ and $\zeta'$  were already in the same free homotopy class in $M_1$, as we wished to show.

From step one, we know that for each $\rho \in i_*(\Theta(\bigwedge(S_1)))$ there is a Reeb orbit $\gamma_\rho$ in $\rho$.

\

\textbf{Step 3:} For each $\rho \in i_*(\Theta(\bigwedge(S_1)))$, the Reeb orbit $\gamma_\rho$ considered in Step 1 is the unique Reeb orbit of $\lambda_{FH}$ in $\rho$.

Let $\gamma$ be a Reeb orbit in $\rho$. If it is contained in $M_1$, we know that $\gamma$ exists also as a closed geodesic in $(S_1,g)$. Using an argument as in step 2 above, it is easy to show that $\gamma$ and $\gamma_\rho$ are freely homotopic in $M_1$, and therefore also in $T_1 S_1$. Projecting to $S_1$ we obtain that $\gamma$ and $\gamma_\rho$ are lifts of geodesics of $(S_1,g)$ in a same free homotopy class of $S_1$. But for each free homotopy class of $S_1$ there is a unique closed geodesic of $(S_1,g)$; this implies that $\gamma= \gamma_\rho$.

Step 3 will now follow if we prove the following claim: \textbf{every Reeb orbit of $\lambda_{FH}$ in $\rho$ is completely contained in $M_1$.}

\textit{Proof of the claim:} if $\gamma$ was contained in $M_2$ then it would be possible to isotopy $\gamma_\rho$ to a curve completely contained in $\partial M_1$. This is impossible by the definition of $\bigwedge(S_1)$.

The only remaining possibility is that $\gamma$ visit both $M_1$ and $M_2$. In this case, it has to visit both $M_1 \setminus N$ and $M_2 \setminus N$ (the reason for that is that if $\gamma$ is completely contained in $M_i \cup N$ convexity of the hyperbolic metric implies that $\gamma$ is in $M_i$).
As $\gamma$ visits both $M_1 \setminus N$ and $M_2 \setminus N$, we know from the Lemma \ref{lemma3} that every curve which is freely homotopic to $\gamma$ has to intersect the torus $\mathbb{T}$. As $\gamma_\rho$ does not intersect $\mathbb{T}$ it cannot be freely homotopic to $\gamma$ which implies that $\gamma \notin \rho$, finishing the proof of step 3.

\

\textbf{Step 4:} End of the proof.

From the previous steps we know that for each $\rho \in i_*(\Theta(\bigwedge(S_1)))$, there exists a unique Reeb orbit $\gamma_\rho \in \rho$. This implies that for each $\rho \in i_*(\Theta(\bigwedge(S_1)))$, the cylindrical contact homology $C\mathbb{H}_{\rho}^{cyl}(M,\lambda_{FH})\neq 0$.

Let $\rho \in i_*(\Theta (\bigwedge(S_1)^{\leq T}))$. Then as we showed, in the previous steps, the unique Reeb orbit of $\lambda_{FH}$ in $\rho$ has action smaller or equal than $T$ and $C\mathbb{H}_{\rho}^{cyl}(M,\lambda_{FH})\neq 0$. This implies that:
\begin{equation}
N^{cyl}_T(\lambda_{FH}) \geq \sharp( i_*(\Theta (\bigwedge(S_1)^{\leq T}))).
\end{equation}

As $i_*$ restricted to $\Theta (\bigwedge(S_1)^{\leq T}))$ is injective, and $\Theta$ is injective we conclude that:
\begin{equation}
\sharp( i_*(\Theta (\bigwedge(S_1)^{\leq T})))= \sharp(\bigwedge(S_1)^{\leq T})\geq e^{aT+b}.
\end{equation}
Combining formulas (77) and (78), we obtain
\begin{equation}
N^{cyl}_T(\lambda_{FH}) \geq e^{aT+b}.
\end{equation}
\qed

\section{Conclusion}\label{section7}

The work of Katok \cite{K} \cite{K2} implies that if an autonomous flow on a 3-manifold has positive topological entropy then there exists a Smale horseshoe as a subsystem of the flow. For a flow a ``horseshoe'' is a compact invariant set where the dynamics are conjugate to that of the suspension of a shift map. In particular, the number of hyperbolic periodic orbits on a ``horseshoe'' of a 3-dimensional flow grows exponentially with respect to the period; see the recent paper \cite{LS} for a refined estimate of this growth. As a consequence, for the contact 3-manifolds $(M,\xi)$ considered in Theorems \ref{theorem3'} and \ref{theorem4'}, we have that for \textbf{every} Reeb flow on $(M,\xi)$ the number of hyperbolic Reeb orbits grows exponentially with the action. This can be summarised by saying that all Reeb flows on these contact manifolds posses a ``complicated'' orbit structure which is forced to exist by the ``complicated'' contact topology of these contact manifolds.

 An interesting property of the entropy estimate used in this paper, \cite{A1} and  \cite{MS} is that it gives estimates on the growth of the number of hyperbolic Reeb orbits also for degenerate contact forms. This kind of information is not obtainable just by studying the growth rate of contact homology.

It is known that the consequences of positivity of topological entropy in higher dimensions are not as strong as in the low dimensional case. In particular positive topological entropy for a flow in dimension bigger than $3$ does not imply the existence of a ``horseshoe'' in the flow. It is however, natural to ask the following question.

\

\textbf{Question 1:} In dimension bigger or equal to $5$, does exponential homotopical growth of periodic orbits for a Reeb flow imply the existence of a compact invariant set where the dynamics are conjugated to a shift?

\

In another direction, one would like to know if it is possible to obtain more dynamical information about the Reeb flows on the contact manifolds covered by Theorems \ref{theorem3'} and \ref{theorem4'}.

\

\textbf{Question 2:} Let $(M,\xi)$ be a manifold satisfying the hypothesis of Theorem \ref{theorem3'} or \ref{theorem4'}, and $\lambda$ a contact form on $(M,\xi)$. Is it true that for the Reeb flow $\phi_{X_\lambda}$ there exists an invariant region of positive measure (with respect to the measure $\lambda\wedge d\lambda$) on which the dynamics of the Reeb flow is ergodic?

\

One important property of many of the contact 3-manifolds covered in Theorem \ref{theorem3'} is that they have positive Giroux torsion. By a theorem of Gay \cite{G} (see also \cite{W}) manifolds with positive Giroux torsion are not strongly fillable. This implies that many of the contact manifolds satisfying the claims of Theorem \ref{theorem3'} are not strongly fillable and therefore different from the unit tangent bundles studied in \cite{MS}, which are exactly fillable. It would be interesting to know if such examples also exist in high dimensions.

\

\textbf{Question 3:} Are there examples of non-symplectically fillable contact manifolds with dimension $\geq 5$, on which every Reeb flow has positive topological entropy? Are there examples in dimensions $\geq 5$, of manifolds which admit infinitely many different contact structures such that on all of them every Reeb flow has positive topological entropy?

\

 We remark also, that in Theorem \ref{theorem3'} we showed the existence of 3-manifolds with hyperbolic components which can be given infinitely many different contact structures whose Reeb flows always have positive topological entropy.  From the perspective of 3-dimensional topology, it would be interesting to have examples of contact structures on hyperbolic 3-manifolds on which every Reeb flow has positive topological entropy.

\

\textbf{Question 4:} Are there examples of contact structures on closed hyperbolic 3-manifolds, on which every Reeb flow has positive topological entropy? Are there hyperbolic 3-manifolds which admit multiple non-diffeomorphic contact structures on which every Reeb flow has positive topological entropy?

\

Lastly we mention that the techniques used in this paper and in \cite{A1}, can also be used in combination with the ideas of Momin \cite{Mo} to establish chaotic behavior of Reeb flows on $(S^3,\xi_{tight})$, when these Reeb flows have a a special link as a Reeb orbit. This and similar results will appear in \cite{A2}.


\begin{thebibliography}{99}


\bibliographystyle{plain}

\bibitem{A3} M. Alves. Growth rate of Legendrian contact homology and dynamics of Reeb flows. \textit{PhD thesis.} Universit\'e Libre de Bruxelles, \textit{December 2014}

\bibitem{A1} M. Alves. Legendrian contact homology and topological entropy. arXiv preprint arXiv:1410.3381 (2014).

\bibitem{A2} M. Alves and P. Salomao. Legendrian contact homology on the complement of Reeb orbits and topological entropy. \textit{preprint in preparation}



\bibitem{B} F. Bourgeois, A survey of contact homology, in ``New perspectives and challenges in symplectic field theory'', 45--71, {\it CRM Proc. Lecture Notes}, {\bf 49} , Amer. Math. Soc., RI, 2009.


\bibitem{CPT} F. Bourgeois, Y. Eliashberg, H. Hofer, K. Wysocki, E. Zehnder.
Compactness results in symplectic field theory. \textit{Geometry and Topology}, 7:799-888 , 2003

\bibitem{BM} F. Bourgeois and K. Mohnke, Coherent orientations in symplectic field theory, \textit{Math. Z.} 248 (2004), no. 1, 123–146

\bibitem{BEE} F. Bourgeois, T. Ekholm,  Y. Eliashberg.  Effect of Legendrian surgery, \textit{Geometry \& topology} (2012) , Volume 16, 1, 301-389

\bibitem{Bo} R. Bowen. Topological entropy and axiom A, \textit{Global Analysis (Proc. Sympos. Pure Math., Vol. XIV, Berkeley, Calif., 1968)}, Amer. Math. Soc., Providence, R.I., 1970, pp. 23–41.


\bibitem{Colin} V. Colin. Une infinit\'e de structures de contact tendues sur les variétés toroidales, \textit{Comentharii Mathematici Elvetici}, 76(2): 353-372, 2001.

\bibitem{CH} V. Colin, K. Honda. Construction control\'ees des champs de Reeb et applications. {\it Geometry and Topology} 9:2193-2226, 2005.

\bibitem{Dr} D. Dragnev, Fredholm theory and transversality for non compact pseudoholomorphic
curves in symplectisations. \textit{Communications on pure and applied mathematics}, 57:726–763,
2004.

\bibitem{SFT} Y. Eliashberg, A. Givental, H. Hofer, Introduction to symplectic field theory,
{\it Geom. Funct. Anal.} (2000), Special Volume, Part II, 560--673.

\bibitem{Fel} A. Fel’shtyn. Dynamical Zeta Functions, Nielsen Theory and Reidemeister torsion,
volume 147 of \textit{Memoirs of the American Mathematical Society}. American Mathematical Society,
Providence, RI, 2000.

\bibitem{Fe} S. Fenley, Homotopic indivisibility of closed orbits of Anosov flows in dimension 3 , \textit{Math. Zeit.} 225 (1997), 289-294

\bibitem{FH} P. Foulon and B. Hasselblatt. Contact Anosov flows on hyperbolic 3–manifolds. \textit{Geometry and Topology} 17 (2013) 1225–1252


\bibitem{FLS} U. Frauenfelder, C. Labrousse, F. Schlenk. Slow volume growth for Reeb flows on spherizations and contact Bott--Samelson theorems. \textit{Preprint} arXiv:1307.7290 (2013).

\bibitem{FS1} U. Frauenfelder, F. Schlenk. Volume growth in the component of the Dehn-Seidel twist. \textit{Geometric \& Functional Analysis GAFA}, v. 15, n. 4, p. 809-838, 2005.

\bibitem{FS2} U. Frauenfelder and F. Schlenk, Fiberwise volume growth via Lagrangian intersections.
\textit{J. Symplectic Geom.} 4 (2006) 117--148.	

\bibitem{FS3} U. Frauenfelder and F. Schlenk., Filtered Hopf algebras and counting geodesic chords. \textit{Mathematische Annalen} 360.3-4 (2014): 995-1020.	


\bibitem{G} D. T. Gay, Four-dimensional symplectic cobordisms containing three-handles, \textit{Geom. Topol.} 10, 1749–1759 (electronic), 2006.

\bibitem{Gr} M. Gromov. Pseudo-holomorphic curves in symplectic manifolds. \textit{Inventiones Mathematicae},
82:307–347, 1985.

\bibitem{HT} M. Handel and W. P. Thurston. Anosov flows on new three manifolds. \textit{Invent. Math.} 59 (1980), no. 2, 95–103.

\bibitem{HK} B. Hasselblatt and A. Katok. Introduction to the Modern Theory of Dynamical Systems (Cambridge, 1995)

\bibitem{H} H. Hofer. Pseudoholomorphic curves in symplectization with applications to the Weinstein
conjecture in dimension three. \textit{Inventiones mathematicae}, 114:515–563, 1993.

\bibitem{P1} H. Hofer, E. Zehnder e K. Wysocki. Properties of pseudoholomorphic curves in symplectizations I: Assymptotics, \textit{Ann. I. H. P. Analyse Non Lin eaire}, 13,(1996),337-379.

\bibitem{P2} H. Hofer, E. Zehnder e K. Wysocki. Properties of pseudoholomorphic curves in symplectizations III: Fredholm theory, \textit{Topics in nonlinear analisys}, Birkhauser, Basel, (1999),381-475.

\bibitem{HMS} U. Hryniewicz, A. Momin, P. A. S. Salom\~ao. A Poincar\'e-Birkhoff theorem for tight Reeb flows on $S^3$, \textit{Inventiones Mathematicae} (2014): 1-90.

\bibitem{J} B. Jiang. Estimation of the number of periodic orbits, \textit{Pacific J. Math.} 172 (1996), 151–
185.

\bibitem{K} A. Katok, Lyapunov exponents, entropy and periodic orbits for diffeomorphisms, \textit{Inst. Hautes Études Sci. Publ. Math.} 51, 137–173, 1980.

\bibitem{K2} A. Katok. Entropy and closed geodesics. \textit{Ergodic Theory Dynam. Systems} 2 339–365 (1983)

\bibitem{LS} Y. Lima, O. Sarig. Symbolic dynamics for three dimensional flows with positive topological entropy, arXiv:1408.3427, 2014.

\bibitem{MS} L. Macarini and F. Schlenk. Positive topological entropy of Reeb flows on spherizations. \textit{Math. Proc. Cambridge Philos. Soc.} 151 (2011) 103--128.

\bibitem{Mo} A. Momin. Contact Homology of Orbit Complements and Implied Existence. \textit{Journal of Modern Dynamics}, Pages: 409 - 472, Issue 3, July 2011

\bibitem{Pat} G. Paternain. {{\it Geodesic Flows}} (Progress in Mathematics) Birkhauser

\bibitem{Rob} C. Robinson, Dynamical Systems (CRC Press, London, 1995)

\bibitem{Vaugon} A. Vaugon. On growth rate and contact homology, \textit{Algebraic Geometry and Topology} 15, no. 2, 623-666, 2015.

\bibitem{W} C. Wendl. Strongly fillable contact manifolds and J-holomorphic foliations, \textit{Duke Math. J.} 151, no. 3, 337–384, 2010.

\end{thebibliography}
\end{document}